\documentclass[11pt]{amsart}

\usepackage[mathscr]{eucal} 
\usepackage{latexsym, amsfonts, amsmath, amssymb, mathrsfs, xspace}
\usepackage[breaklinks=true]{hyperref}

\def\k{{\Bbbk}}

\newcommand{\gr}{\mathrm{gr}}
\newcommand{\lex}{\mathrm{lex}}
\newcommand{\calo}{\ensuremath{\mathcal{O}}\xspace}
\newcommand{\mfsl}{\mathfrak{sl}}
\newcommand{\inv}{\mathsf{i}}
\newcommand{\eps}{\epsilon}
\newcommand{\Z}{\mathbb{Z}}
\newcommand{\zz}{\Z_{\ge 0}}
\newcommand{\mf}[1]{\ensuremath{\mathfrak{#1}}\xspace}
\newcommand{\brac}[1]{\ensuremath{\{#1\}}\xspace}
\newcommand{\uqsl}{\ensuremath{U_q(\mf{sl}_2)}\xspace}
\newcommand{\usl}{\ensuremath{{U}(\mf{sl}_2)}\xspace}
\newcommand{\onto}{\twoheadrightarrow}

\newcommand{\kk}{{\mathit{k}}}
\newcommand{\oH}{{\bar{H}}}
\newcommand{\oX}{{\bar{X}}}
\newcommand{\oY}{{\bar{Y}}}
\newcommand{\oE}{{\bar{E}}}
\newcommand{\oF}{{\bar{F}}}
\newcommand{\oK}{{\bar{K}}}

\newcommand{\cz}{\ensuremath{\eta_0(C_0)}\xspace}
\newcommand{\cm}{\ensuremath{\eta_{-1}(C_0)}\xspace}
\newcommand{\bz}{\ensuremath{\eta_0(b_1)}\xspace}
\newcommand{\bm}{\ensuremath{\eta_{-1}(b_1)}\xspace}
\newcommand{\az}{\ensuremath{\eta_0(a)}\xspace}
\newcommand{\am}{\ensuremath{\eta_{-1}(a)}\xspace}
\newcommand{\amm}{\ensuremath{\eta_{-2}(a)}\xspace}

\newcommand{\ddeg}{\mathrm{deg}}

\newcommand{\into}{\hookrightarrow}
\newcommand{\iso}{\stackrel{\sim}{\longrightarrow}}
\newcommand{\mapdef}[1]{\ensuremath{\overset{#1}{\longrightarrow}}\xspace}

\newcommand{\mcA}{\mathcal{A}}
\newcommand{\F}{\mathscr{F}}

\DeclareMathOperator{\Ext}{\ensuremath{Ext^1_\calo}}
\DeclareMathOperator{\Hom}{\ensuremath{Hom}}
\DeclareMathOperator{\im}{im}
\newcommand{\qdiff}{(q^{-2} + q^2)}
\newcommand{\hhom}{\ensuremath{\Hom_\calo}\xspace}

\theoremstyle{plain}
\newtheorem{definition}[equation]{Definition}
\newtheorem{corollary}[equation]{Corollary}
\newtheorem{lemma}[equation]{Lemma}
\newtheorem{proposition}[equation]{Proposition}
\newtheorem{theorem}[equation]{Theorem}

\theoremstyle{definition}
\newtheorem{remark}[equation]{Remark}
\newtheorem*{acknowledgement}{Acknowledgements}
\newtheorem*{stand}{Standing Assumption}

\numberwithin{equation}{section}
 
\begin{document}

\title{Quantized symplectic oscillator algebras of rank one}

\author{Wee Liang Gan}
\address{Department of Mathematics, Massachusetts Institute of
Technology, Cambridge, MA 02139, U.S.A.}
\email{wlgan@math.mit.edu}

\author{Nicolas Guay}
\address{Department of Mathematics, University of Illinois at
Urbana-Champaign, Urbana, IL 61801, U.S.A.}
\email{nguay@math.uiuc.edu}

\author{Apoorva Khare}
\address{Department of Mathematics, University of Chicago,
Chicago, IL 60637, U.S.A.}
\email{apoorva@math.uchicago.edu}



%


\begin{abstract}
A quantized symplectic oscillator algebra of rank 1 is a
PBW deformation of the smash product of the quantum plane
with \uqsl.
We study its representation theory, and in particular,
its category \calo. 
\end{abstract}
\maketitle

\section{Introduction}
\label{intro}

Let $V$ be a finite dimensional complex vector space
equipped with a nondegenerate skew-symmetric bilinear form.
In \cite[\S 4]{EGG}, Etingof, Gan and Ginzburg introduced the
family of infinitesimal Hecke algebras $\mathscr H_\beta$
associated to $\mathfrak{sp}(V)$. The algebras
$\mathscr H_\beta$ are PBW deformations of 
$\mathbb C[V]\rtimes U(\mathfrak{sp}(V))$. 
On one hand, they are similar to the symplectic reflection
algebras introduced by Etingof and Ginzburg in \cite{EG}
(and by Crawley-Boevey and Holland in 
\cite{CBH} when $\dim V$ is 2). 
On the other hand, they are also similar to universal
enveloping algebras of Lie algebras.
In the case when $\dim V$ is 2, the algebra
$\mathscr H_\beta$ was also called a symplectic oscillator
algebra in \cite{Kh} (see \cite[Example 4.12]{EGG}); 
we shall refer to the $\mathscr H_\beta$
in this case as the symplectic oscillator algebras of rank 1.

The representation theory of the symplectic oscillator algebras of rank 1
was studied by Khare in \cite{Kh}. 
In our present paper, we show that the main results of \cite{Kh} can
naturally be $q$-deformed.
One of our main results is that, in the $q$-deformed setting, there exist
PBW deformations whose finite dimensional representations are completely
reducible. 
(The same proof can also be adapted to the original setting in
\cite{Kh}.)

Fix a ground field $\k$, with char $\k \neq 2$, 
and an element $q\in\k^\times$ such that $q^2\neq 1$.
Since the quantum plane $\k_q[X,Y]:=\k\langle X,Y\rangle/(XY-qYX)$
is a module-algebra over the Hopf algebra
$\uqsl$, one can form the smash product algebra
$\k_q[X,Y]\rtimes \uqsl$, cf. \cite{Mo}.
Our main object of study is a deformation of this algebra,
defined for each element $C_0$ in the center of $\uqsl$ as follows.

\begin{definition} \label{maindef}
The {\it quantized symplectic oscillator algebra of rank $1$} is the
algebra $A$ generated over $\k$ by the elements
$E, F, K, K^{-1}, X, Y$ with defining relations
\begin{gather}
KEK^{-1} = q^2 E, \quad KFK^{-1} = q^{-2} F, \quad
[E,F] = \frac{K-K^{-1}}{q-q^{-1}}, \label{r1}\\
EX = qXE, \quad EY = X + q^{-1}YE, \label{r2}\\
FX = YK^{-1} + XF, \quad FY= YF, \label{r3}\\
KXK^{-1} = qX, \quad KYK^{-1} = q^{-1}Y, \label{r4}\\
qYX-XY = C_0. \label{r5}
\end{gather}
\end{definition}

The PBW Theorem for $A$ is the statement that
the set of elements $F^aY^bK^cX^dE^e$
(for $a,b,d,e \in \Z_{\ge 0}$, $c\in \Z$) form
a basis for $A$.
We will prove this in \S \ref{sPBW}.
Let us make some comments on Definition \ref{maindef}.

\smallskip
\begin{remark}\hfill
\begin{enumerate}
\item Observe that the subalgebra of $A$ generated by 
$E, F, K$ and $K^{-1}$ is isomorphic to $\uqsl$.
When $C_0=0$, the algebra $A$ is $\k_q[X,Y]\rtimes \uqsl$.

\item In \cite{Kh}, the (deformed) symplectic oscillator algebra $H_f$ is
defined, for each polynomial $f\in \k[t]$, to be the quotient of
$T(V) \rtimes U(\mathfrak{sp}(V))$ by the relations
$[y,x] = \omega(x,y)(1+f(\Delta))$ for all $x,y\in V$,
where $\Delta$ is the Casimir element in $U(\mathfrak{sp}(V))$.
The PBW Theorem for $H_f$ was proved in 
\cite[Theorem 9]{Kh} when $\dim V=2$.
However, it is not true in general when $\dim V>2$.
The formula for obtaining PBW deformations of $\mathbb C[V]\rtimes
U(\mathfrak{sp}(V))$ is given in \cite[Theorem 4.2]{EGG}.

For the rest of this paper, $H_f$ will always mean the case $\dim V=2$.

\item The symplectic oscillator algebra when $\dim V=2$ is analogous to
the algebra 
$$\frac{\k\langle X,Y\rangle \rtimes \k[\Gamma]}{(YX-XY-\zeta)}$$

\noindent where $\Gamma$ is a finite subgroup of $SL(V)$ and $\zeta$ is
an element in the center of $\k[\Gamma]$, introduced and studied by
Crawley-Boevey and Holland in \cite{CBH}.    \hfill$\Diamond$
\end{enumerate}
\end{remark}

The algebra $A$ is very similar to quantized universal enveloping
algebras of semisimple Lie algebras in many ways. For example, we can
construct Verma modules using the PBW Theorem, 
define highest weight modules, and study its category \calo. 
The main results of the paper are the following:
\begin{itemize}
\item Necessary and sufficient conditions for a simple highest
weight module to be finite dimensional. (Theorem \ref{fd})

\item A description of the Verma modules of $A$ using 
the Verma modules of its subalgebra $U_q(\mathfrak{sl}_2)$.
(Theorem \ref{nonint} and Theorem \ref{int})

\item A block decomposition for \calo, and
a proof that \calo is a highest weight category.
(Corollary \ref{blockd} and Proposition \ref{hwc})

\item Necessary and sufficient conditions for the finite
dimensional representations to be completely reducible.
(Theorem \ref{weyl})

\item A proof that the center of $A$ is trivial when
$C_0\neq 0$. (Theorem \ref{thmcenter})
\end{itemize}

Note that since the center of $A$ is trivial when
$C_0\neq 0$, the original approach in \cite{BGG1} for
the decomposition of \calo does not work for our algebra.

\smallskip

\noindent {\bf Organization of the paper.}
We prove the PBW Theorem in \S \ref{sPBW}.
Then we study in \S \ref{scyclicmod} the actions of the ``raising''
operators $E$ and $X$ on highest weight modules.
Thenceforth we assume that $q$ is not a root-of-unity.
In \S \ref{sfinite}, we determine necessary and sufficient conditions for
a simple highest weight module to be finite dimensional.
In \S \ref{ssubquot}, we determine conditions for existence of maximal
vectors in Verma modules.
Beyond this point, we assume that $C_0 \neq 0$.
Then we study those Verma modules in \S \ref{snonint} whose highest
weights are not of the form $\pm q^n$, where $n\in\Z$.
In \S \ref{sint} we study Verma modules whose highest weights are of the
form $\pm q^n$ where $n\in\Z$. We obtain, in the following section, a
decomposition of the category \calo into blocks, each of which is a
highest weight category.
In \S \ref{sblockdec}, we show that various ways of decomposing \calo
into blocks are actually equivalent.
A characterization of all cases when complete reducibility holds is the
content of \S \ref{sweyl}.
The proof of the complete reducibility in this section makes use of
our results obtained in the earlier sections, in particular,
the decomposition of \calo.
In \S \ref{scenter}, we prove that the center of $A$ is trivial if
$C_0\neq 0$. The next section contains some more  
results about the Verma modules.
Finally, we explain how to take the classical limit $q\to 1$ to obtain
the algebra $H_f$ and its highest weight modules in \S \ref{sclaslimit}.

\section{PBW Theorem} \label{sPBW}

The relations (\ref{r2}), (\ref{r3}) and (\ref{r4})
imply that $qYX-XY$ commutes with $E, F, K$ and $K^{-1}$.
However, $C_0$ does not necessarily commute with $X, Y$.

\begin{theorem}\label{PBW}
The set of elements $F^aY^bK^cX^dE^e$,
where $a,b,d,e \in \Z_{\ge 0}$, $c\in \Z$, is a basis for $A$.
\end{theorem}

\begin{proof}
We shall use the Diamond Lemma; see \cite[Theorem 1.2]{Be} or \cite{BG}.

To be precise, we write $K^{-1}$ as $L$, so 
\begin{equation} \label{r6}
KL=LK=1.
\end{equation}

We now define a semigroup partial ordering $\le$ on the set $W$ 
of words in the generators $E,F,K,L,X,Y$. First, define the
lexicographic ordering $\le_{\lex}$ on $W$ by ordering the generators 
as $F,Y,L,K,X,E$. For each word $w\in W$, let $n(w)$ be the total
number of times $X$ and $Y$ appear in $w$.
Now, given two words $w$ and $u$, we define $w\le u$ if
\begin{itemize}
\item $n(w)<n(u)$, or
\item $n(w)=n(u)$ and ${\rm length}(w)<{\rm length}(u)$, or
\item $n(w)=n(u)$, ${\rm length}(w)={\rm length}(u)$ and
$w\le_{\lex}u$.
\end{itemize}

This is a semigroup partial ordering 
which satisfies the descending chain condition 
and is also compatible with the reduction 
system given by our relations (\ref{r1})--(\ref{r5}) 
and (\ref{r6}). We have to check that the 
ambiguities are resolvable, which we do below.
The Diamond Lemma then implies that the irreducible words
\begin{gather*}
\{F^aY^bL^cX^dE^e\mid a,b,d,e\in\Z_{\ge 0}, c\in\Z_{>0}\} \\
\bigcup \{F^aY^bK^cX^dE^e\mid a,b,c,d,e\in\Z_{\ge 0}\}
\end{gather*}
form a basis for $A$.

Here are the details of the verification. Let us first 
write down our reduction system:
$$EK \to q^{-2}KE,\quad KF\to q^{-2}FK,\quad LK\to 1,
\quad KL\to 1,$$
$$EF\to FE+(K-L)/(q-q^{-1}),\quad EX\to qXE,\quad
EY\to X+q^{-1}YE,$$
$$ XF\to FX-YL,\quad YF\to FY, \quad XY\to qYX-C_0,$$
$$ EL\to q^2LE,\quad LF\to q^2 FL,\quad XK\to q^{-1}KX,$$
$$ KY \to q^{-1}YK,\quad XL\to qLX,\quad LY\to qYL.$$

Observe that there is no inclusion ambiguity and all
overlap ambiguities appear in words of length 3.
Moreover, if $X$ and $Y$ do not appear in a word, 
then it is reduction unique by the PBW theorem for $U_q(\mfsl_2)$.
Thus, the words which we have to check are:
$$LYF, KYF, XYF, EYF, EXF, XLF, XKF, KLY, $$
$$XLY, ELY, XKY, EKY, EXY, XKL, EXL, EXK. $$

We now show that all these ambiguities are resolvable:
$$L(YF) \to (LF)Y \to q^2F(LY) \to q^3FYL$$
$$(LY)F \to qY(LF) \to q^3(YF)L \to q^3FYL$$

$$K(YF) \to (KF)Y \to q^{-2}F(KY) \to q^{-3}FYK$$
$$(KY)F \to q^{-1}Y(KF) \to q^{-3}(YF)K \to q^{-3}FYK$$

$$X(YF) \to (XF)Y \to F(XY)-Y(LY) \to qFYX-FC_0-qYYL$$
$$XYF \to qYXF-C_0F \to qYFX-qYYL-C_0F \to qFYX-qYYL-C_0F$$

$$E(YF) \to (EF)Y \to F(EY)+(KY-LY)/(q-q^{-1})$$
$$\to FX+q^{-1}FYE+(q^{-1}YK-qYL)/(q-q^{-1})$$
$$(EY)F \to XF+q^{-1}Y(EF) $$
$$\to FX-YL+q^{-1}(YF)E+(q^{-1}YK-q^{-1}YL)/(q-q^{-1}) $$
$$\to FX-YL+q^{-1}FYE+(q^{-1}YK-q^{-1}YL)/(q-q^{-1})$$

$$E(XF) \to (EF)X-(EY)L $$
$$\to F(EX)+(KX-LX)/(q-q^{-1})-XL-q^{-1}Y(EL)$$
$$\to qFXE+(KX-LX)/(q-q^{-1}) -qLX -qYLE$$
$$(EX)F \to qX(EF) \to q(XF)E+(qXK-qXL)/(q-q^{-1})$$
$$\to qFXE-qYLE+(KX-q^2LX)/(q-q^{-1})$$

$$X(LF) \to q^2(XF)L \to q^2F(XL)-q^2YLL \to q^3FLX-q^2YLL$$
$$(XL)F \to qL(XF) \to q(LF)X-q(LY)L \to q^3FLX -q^2YLL$$

$$X(KF) \to q^{-2}(XF)K \to q^{-2}F(XK)-q^{-2}YLK \to
q^{-3}FKX-q^{-2}Y(LK)$$
$$XKF\to q^{-1}KXF \to q^{-1}KFX-q^{-1}KYL
\to q^{-3}FKX-q^{-2}YKL$$

$$K(LY) \to q(KY)L \to Y(KL)$$
$$(KL)Y \to Y$$

$$X(LY) \to q(XY)L \to q^2Y(XL)-qC_0L \to q^3YLX-qC_0L$$
$$(XL)Y \to qL(XY) \to q^2(LY)X-qLC_0 \to q^3YLX-qLC_0$$

$$E(LY) \to q(EY)L \to qXL+Y(EL) \to q^2LX+q^2YLE$$
$$(EL)Y \to q^2L(EY) \to q^2LX+q(LY)E \to q^2LX+q^2YLE$$

$$X(KY) \to q^{-1}(XY)K \to Y(XK)-q^{-1}C_0K
\to q^{-1}YKX-q^{-1}C_0K$$
$$(XK)Y \to q^{-1}K(XY) \to (KY)X-q^{-1}KC_0
\to q^{-1}YKX-q^{-1}KC_0$$

$$E(KY) \to q^{-1}(EY)K \to q^{-1}XK+q^{-2}Y(EK)
\to q^{-2}KX+q^{-4}YKE$$
$$(EK)Y \to q^{-2}K(EY) \to q^{-2}KX+q^{-3}(KY)E
\to q^{-2}KX+q^{-4}YKE$$

$$E(XY) \to q(EY)X-EC_0 \to qXX+Y(EX)-EC_0
\to qXX+qYXE-EC_0$$
$$(EX)Y \to qX(EY) \to qXX+(XY)E \to qXX+qYXE-C_0E$$

$$X(KL) \to X$$
$$(XK)L \to q^{-1}K(XL) \to (KL)X \to X$$

$$E(XL) \to q(EL)X \to q^3L(EX) \to q^4LXE$$
$$(EX)L \to qX(EL) \to q^3(XL)E \to q^4LXE$$

$$E(XK) \to q^{-1}(EK)X \to q^{-3}K(EX) \to q^{-2}KXE$$
$$(EX)K \to qX(EK) \to q^{-1}(XK)E \to q^{-2}KXE$$
This completes the proof of Theorem \ref{PBW}. 
\end{proof}

\noindent This method can also be applied to $H_f$, and provides a
simpler proof than in \cite{Kh}.

We may define a $\Z_{\ge 0}$-filtration on $A$ by assigning
$\ddeg E=\ddeg F=1$, $\ddeg K=\ddeg K^{-1}=0$, 
and $\ddeg X=\ddeg Y$ to be some sufficiently big number
so that, by Theorem \ref{PBW}, the associated graded algebra
$\gr A$ is a skew-Laurent extension
of a quantum affine space, cf. e.g. \cite{BG}.
Hence, we obtain the following corollary.

\begin{corollary}\label{noeth} 
The algebra $A$ is a Noetherian domain.
\end{corollary}

\section{Standard cyclic modules} \label{scyclicmod}

Given a $\k[K,K^{-1}]$-module $M$ and $a \in \k^\times$, we define
$M_a = \{ m \in M : K \cdot m = am \}$ and denote by $\Pi(M)$ the set of
{\it weights}: $\{a \in \k^\times : M_a \neq 0 \}$.
We consider $A$ as a $\k[K,K^{-1}]$-module on which $K^c$ ($c\in\Z$) acts
by conjugation.

\begin{lemma}\hfill
\begin{enumerate}
\item If $M$ is a $\k[K,K^{-1}]$-module,
then the sum $\sum_{a \in \k^\times} M_a$ is direct, and $K$-stable.

\item If $M$ is any $A$-module, then $A_a M_b \subset M_{ab}$.

\item We have: $A = \bigoplus_a A_a$, and $\k[K,K^{-1}] \subset A_1$.
\end{enumerate}
\end{lemma}

Note that $A$ contains subalgebras
$B_+ = \langle E,X,K,K^{-1} \rangle$ 
and $B_- = \langle F,Y,K,K^{-1} \rangle$. 
We define $N_+$ (resp. $N_-$) to be the nonunital subalgebra of $A$
generated by $E,X$ (resp. $F,Y$). These are analogs of the enveloping
algebras of Borel or nilpotent subalgebras of a semisimple Lie algebra.

Later on, we will use often the ``purely CSA'' (CSA stands for Cartan
subalgebra) map $\xi : A \to \k[K,K^{-1}]$ defined as follows: write each
element $U \in A$ in the PBW basis given in Theorem \ref{PBW}, then
$\xi(U)$ is the sum of all vectors in $\k [K,K^{-1}]$, i.e. 
$U - \xi(U) \in N_-A + AN_+$.

We need some terminology that is standard in representation theory. If
$M$ is an $A$-module, a {\it maximal vector} is any $m \in M$ that is
killed by $E,X$ and is an eigenvector for $K,K^{-1}$. A {\it standard
cyclic module} is one that is generated by exactly one maximal vector.
For each $r\in\k^\times$, define the {\it Verma module} $Z(r) := A /(AN_+
+ A(K - r \cdot 1))$, cf. \cite{H,Kh}. It is a free $B_-$-module of rank
one, by the PBW theorem for $A$, hence isomorphic to $\k[Y,F]$ and has a
basis $\{F^iY^j : i,j \geq 0\}$. Furthermore, $\Pi(Z(r)) = \{ q^{-n}r,\ n
\geq 0\}$. 

The proof of the following proposition is standard --
see e.g. \cite{H} or \cite{Kh}.

\begin{proposition}\hfill
\begin{enumerate}
\item $Z(r)$ has a unique maximal submodule $W(r)$,
and the quotient $Z(r)/W(r)$ is a simple module $V(r)$.
\item Any standard cyclic module is a quotient of some Verma module.
\end{enumerate}
\end{proposition}

We may identify $\uqsl$ with the subalgebra of $A$
generated by $E, F, K$ and $K^{-1}$.
Let $\mf{Z}(\uqsl)$ denote the center of \uqsl, and denote by $Z_C(r)$
and $V_C(r)$ the Verma and simple \uqsl-module respectively, of highest
weight $r \in\k^\times$. We note that any $z \in \mf{Z}(\uqsl)$ acts on
any standard cyclic \uqsl-module with highest weight $r$ by the scalar
$\xi(z)(r)$, where we evaluate the (finite) Laurent polynomial $\xi(z)
\in \k [K,K^{-1}]$ at $r \in \k^\times$. Define $c_{0r} = \xi(C_0)(r)$ to
be the scalar by which $C_0$ acts on a \uqsl-Verma module $Z_C(r)$.

Now we introduce some more notation. We know that units in the
$\k$-algebra $\k[K,K^{-1}]$ are all of the form $b K^m$, where $b \in
\k^\times$ and $m \in \Z$. We denote this set (of all units) by
$\k^{\times} K^\Z$. Moreover, for any $a\in\k^\times K^\Z $, define
$\brac{a} := \frac{a - a^{-1}}{q - q^{-1}}$. The following identity is
now easy to check:
\begin{equation}\label{s1}
a \brac{b} - b\brac{a} = \brac{a^{-1}b}\quad \mbox{ for all } a, b \in
\k^\times K^\Z.
\end{equation}

We use the identity in proving the next result, as well as Theorem
\ref{s9} below.

\begin{lemma}\label{s2}
Suppose $a,b \in \k^\times K^\Z$. Then 
\begin{itemize}
\item[(a)] $\brac{a^{-1}} = -\brac{a}$, and 
\item[(b)] $q^{-1} \brac{b} + b = \brac{qb}$.
\end{itemize}
\end{lemma}

\begin{proof}
To prove (a), we set $b=1$ in \eqref{s1}, 
and to prove (b), we set $a = q^{-1}$ in \eqref{s1}.
\end{proof}

We shall write $Z(r) \to V \to 0$ to
mean that $V$ is a standard cyclic $A$-module 
with highest weight $r$.

As we shall see, many standard cyclic (resp. Verma, simple) $A$-modules
$Z(r) \to V \to 0$ are a direct sum of a progression of standard cyclic
(resp. Verma, simple) \uqsl-modules of highest weight $t = r, q^{-1}r,
\dots$, each such module having multiplicity one as well. The specific
equations governing such a direct sum $V = \oplus_i V_{C,q^{-i}r}$ are
the subject of this section.

For $m\ge 2$, we define
\begin{equation}\label{s3}
\alpha_{r,m} = \sum_{j=0}^{m-2} \brac{q^{1-j}r}c_{0,q^{-j}r} .
\end{equation}

\noindent This constant will play a fundamental role in the 
rest of this paper.
(We remark that this constant $\alpha_{rm}$ is different 
from the constant that was also
 denoted by $\alpha_{rm}$ in \cite{Kh}.)

Let $\eps = \pm 1$ henceforth. 
We will also need the constant 
$$d_{r,m} := \frac{\alpha_{r,m}}{\brac{q^{2-m}r}\brac{q^{3-m}r}},$$ 

\noindent which is defined for all $m$, if $r$ is not of the form $\eps
q^l$, or for $2\le m\le l+1$, if $r= \eps q^l$ (where $l\in\Z$).

\begin{lemma}\label{s7}
Given $r \in \k^\times$ and $n \in \zz$, whenever all terms below are
defined, we have
$$\brac{q^{1-n}r} d_{r,n+1} = \brac{q^{3-n}r} d_{r,n} + c_{0,q^{1-n}r}.$$
\end{lemma}

\begin{proof}
We have
\begin{align*}
\brac{q^{2-n}r}(\brac{q^{1-n}r} d_{r,n+1}) =& \alpha_{r,n+1} \\
=& \alpha_{r,n} + \brac{q^{2-n}r}c_{0,q^{1-n}r} \\
=& \brac{q^{2-n}r}(\brac{q^{3-n}r}d_{r,n} + c_{0,q^{1-n}r})
\end{align*}
Since all terms in the claim are defined, $\brac{q^{2-n}r} \neq 0$ and
can be cancelled from both sides.
\end{proof}

\noindent We now imitate the structure theory in \cite[\S 9]{Kh}.

\begin{theorem}\label{s5}
Let $V = A v_r$ be a standard cyclic module, where $v_r$ is a highest
weight vector of weight $r \in \k^\times$. Suppose that $r\neq q^j$ for
$1\le j\le m-1$, and where $m \in \zz$. 
Then we have the following:
\begin{enumerate}
\item $v_r$ and $v_{q^{-1}r} := Yv_r$ are \uqsl-maximal vectors.

\item Suppose $1\le n\le m$. Set $t_n= q^{-n}r$.  
Define inductively:

\begin{equation}\label{Sn}
v_{t_n} := 
Yv_{t_{n-1}} + d_{r,n}
Fv_{t_{n-2}} 
\end{equation}
If $n\ge 2$, the following two equalities hold:
\begin{equation}\label{Rn}
Xv_{t_{n-1}} = EY v_{t_{n-1}} = -\frac{\alpha_{r,n}}{\brac{t_{n-3}}}
v_{t_{n-2}} 
\end{equation}
Moreover, $v_{t_n}$ is \uqsl-maximal, i.e. $Ev_{t_n}=0$. It is a maximal
vector for the algebra $A$ if and only if
$\alpha_{r,n+1}=0$.
\item There exist \emph{monic} polynomials 
$$p_{r,n}(Y,F) = Y^n + c_1 F Y^{n-2} + c_2 F^2 Y^{n-4} + \dots \quad 
(\mbox{where } c_i \in \k)$$ 

\noindent that satisfy $p_{r,n}(Y,F) v_r = v_{t_n}$.
\end{enumerate}
\end{theorem}

\begin{proof}
The last part is obvious from the defining equations, so we show the rest
now.
\begin{enumerate}
\item $v_r$ is $A$-maximal and hence $Ev_r = 0$. Similarly, $EYv_r = Xv_r
= 0$.

\item We proceed by induction, so we assume that all the statements are
true when $n=k$ and we want to show that they are true when $n=k+1$.

\noindent (a) By induction, $v_{t_{k}}$ is \uqsl-maximal, so $Xv_{t_k} =
(EY-q^{-1}YE)v_{t_k} = EYv_{t_k}$.

\noindent (b) If $n$ is $0$ or $1$, 
then we are done from the first part (since we may
choose to set $v_{t_{-1}} = 0$ if we wish). 
If $n=k+1$ and $k>1$, we have
\begin{align*}
Xv_{t_k} = & X(Yv_{t_{k-1}} + d_{r,k}Fv_{t_{k-2}}) \\
 = & (qYX-C_0)v_{t_{k-1}} + d_{r,k}(FX-YK^{-1})v_{t_{k-2}}
\end{align*}

Using the induction hypothesis, we get
\begin{align*}
Xv_{t_k}
 = & qY(-d_{r,k}\brac{t_{k-2}})v_{t_{k-2}} - c_{0,t_{k-1}}v_{t_{k-1}} \\
 & + d_{r,k}(-Fd_{r,k-1}\brac{t_{k-3}}v_{t_{k-3}} -
	Y(t_{k-2})^{-1}v_{t_{k-2}})
\end{align*}

Regrouping terms, we then have
\begin{align*}
Xv_{t_k}
 = & -d_{r,k}Yv_{t_{k-2}}(q \brac{t_{k-2}} + (t_{k-2})^{-1})\\
& -  d_{r,k}\brac{t_{k-3}}(d_{r,k-1}Fv_{t_{k-3}}) -
	c_{0,t_{k-1}}v_{t_{k-1}}
\end{align*}

Now use Lemma \ref{s2} and regroup terms to get
$$Xv_{t_k}
= -d_{r,k} \brac{t_{k-3}}(Yv_{t_{k-2}} + d_{r,k-1}Fv_{t_{k-3}}) -
	c_{0,t_{k-1}} v_{t_{k-1}}$$

Applying the induction hypothesis again, we get
$$Xv_{t_k}
\displaystyle = -d_{r,k} \brac{t_{k-3}}v_{t_{k-1}} -
	c_{0,t_{k-1}}v_{t_{k-1}} =
	-\frac{\alpha_{r,k+1}}{\brac{t_{k-2}}}v_{t_{k-1}}$$

The last equality here uses equation \eqref{s3} and Lemma \ref{s7}.
This completes the induction.

\noindent (c) By induction, $v_{t_k}$ is killed by $E$, so 
\[ EYv_{t_k} = q^{-1}Y(Ev_{t_k}) + Xv_{t_k} = Xv_{t_k} = -d_{r,k+1}
	\brac{t_{k-1}} v_{t_{k-1}}\]

\noindent and 
$$EFv_{t_{k-1}} = (FE + \brac{K})v_{t_{k-1}} =
\brac{t_{k-1}}v_{t_{k-1}}.$$ 
Hence, the vector $Yv_{t_k} +
d_{r,k+1}Fv_{t_{k-1}}$ is indeed killed by $E$. In other words,
$v_{t_{k+1}}$ is a maximal \uqsl-vector.

\noindent Finally, $v_{t_k}$ is $A$-maximal if and only if
$Xv_{t_k} =0$,  which holds
if and only if $\alpha_{r,k+1} = 0$
(using
\eqref{Rn} for $n=k+1$ and note that $\brac{t_{k-2}} \neq 0$).
\end{enumerate}
\end{proof}

\noindent{\bf Example.}
Let us take a look at the undeformed case $C_0=0$. The following proposition
holds under this assumption.

\begin{proposition} Assume $C_0=0$.
\begin{enumerate}
\item Every Verma module $Z(r)$ is a direct sum $Z(r) = \bigoplus_{n \geq
0} Z_C(q^{-n}r)$ of \uqsl-Verma modules. It has a submodule $Z(q^{-1}r)$,
and the quotient $Z_C(r)$ is annihilated by $X,Y$.
\item The simple module $V(r)$ equals $V_C(r)$ and is annihilated by
$X,Y$.
\end{enumerate}
\end{proposition}

\begin{proof}\hfill
\begin{enumerate}
\item We claim that the structure equations, analogous to those in
Theorem \ref{s5}, now become
$$v_{q^{-n}r} = Y^n v_r; \qquad Xv_{q^{-n}r} = Ev_{q^{-n}r} = 0$$
Firstly, $X$ commutes with $Y$ since $C_0=0$, so 
we have $$X(Y^n v_r)
= Y^n (Xv_r) = 0.$$ 
Next, $$EY^n v_r = q^{-1} Y (EY^{n-1} v_r) + XY^{n-1}
v_r = 0$$ by induction on $n$.
Hence, each $Y^n v_r$ is maximal.
We have $$Z(q^{-1}r)\iso A\cdot Yv_r\into Z(r).$$
We also
have isomorphisms $$Z(r) / Z(q^{-1}r) \cong \sum_{n \geq 0} \k F^n
v_r \cong Z_C(r).$$

Now $YF^nv_r=F^nYv_r \in Z(q^{-1}r)$, hence $YF^n\bar{v}_r=0$, where
$\bar{v}_r$ is the image of $v_r$ in the quotient $Z(r)/Z(q^{-1}r)$.
Moreover, 
$$XF^n\bar{v}_r=FXF^{n-1}\bar{v}_r-YK^{-1}F^{n-1}\bar{v}_r=0$$ 
by induction. This proves the last claim of part (1).

\item Since $Z(r)/Z(q^{-1}r)$ is annihilated by $X$ and $Y$, the maximal
submodule of $Z(r)$ corresponds, in this quotient, to the maximal
\uqsl-submodule of $Z_C(r)$.  
\end{enumerate}
\end{proof}

%
%

\begin{stand}
From now on, unless otherwise stated, we assume that $q$ is {\it not} a
root of unity.  \hfill$\Diamond$
\end{stand}

\smallskip

In this case, the center $\mathfrak{Z}(\uqsl)$ is generated by the
Casimir element
$$ C := FE + \frac{qK+q^{-1}K^{-1}}{(q-q^{-1})^2} $$
and we will write $C_0 = p(C)$ for some polynomial $p\in\k[t]$.
We also let $c_r = \xi(C)(r)$ to be the scalar by which $C$ acts
on the \uqsl-Verma module $Z_C(r)$.
Thus, 
\begin{equation} \label{eqcocor}
c_r = (qr + q^{-1}r^{-1}) / (q - q^{-1})^2
\qquad \mbox{ and }\qquad c_{0r} = p(c_r).
\end{equation}

The following proposition will play an important role 
in obtaining the
decomposition of category
\calo (which will be defined later)
and in proving that Verma modules have finite
length.

\begin{proposition}\label{polynom1}
If $q$ is not a root of unity, then
$\alpha_{r,m}$ is of the form
$(q^mr)^{-N}b(r,q^m)$ for some polynomial $b \in \k[S,T]$, 
and some $N \in \Z_{>0}$.
\end{proposition}

\begin{proof}
It is clear from (\ref{eqcocor})
that $c_{0,r}= h(qr)$ for
some $h \in \k[T,T^{-1}]$. Now $g(T) := h(T) \cdot \brac{T}$ is
also in $\k[T,T^{-1}]$. We observe, from equation \eqref{s3}, that
$$\alpha_{r,m} = \sum_{j=0}^{m-2}g(q^{1-j}r).$$
We write $h(T) = \sum_{i=-M}^M b_i T^i$.
By \cite[Lemma 2.17]{Ja},
we have $b_i = b_{-i}$ for each $i$. 
Hence by definition of $g$, if $g(T) = \sum_{i=-N}^N a_i T^i$ (where
$N=M+1$), then $a_{-i} = -a_i$. In particular, $a_0 = 0$.

Recall that we are assuming that $q$ is not a root of unity.
Interchanging the finite sums, we get
\begin{eqnarray*}
\alpha_{r,m} & = & \sum_{j=0}^{m-2} \sum_{i=-N}^N a_i r^i q^{i(1-j)} =
\sum_{i=-N}^N a_i r^i \sum_{j=0}^{m-2} (q^{-i})^{j-1} \\
 & = & \sum_{i=-N}^N \frac{a_i}{q^{-i}-1} r^i (q^{i-mi}-1)q^i
\end{eqnarray*}
\noindent Henceforth, denote by $\sum'$ the summation above
with the $i=0$ term omitted.
If we set $b(S,T) = 
\underset{i=-N}{\overset{N}{\sum'}}\ \displaystyle 
\frac{a_i}{q^{-i}-1}q^i S^{N+i}(q^i T^{N-i} - T^N)$, then $\alpha_{r,m} =
(q^m r)^{-N}b(r,q^m)$.
Here, $b$ is a polynomial in $S,T$, and we are done.
\end{proof}

\section{Finite dimensional modules}\label{sfinite}

We will first give an example for which the category of finite
dimensional modules over $A$ is not semisimple. Afterwards, assuming that
$q$ is not a root of unity, we will give a (rough) classification of all
simple finite dimensional ($K$-semisimple) modules.

\medskip

\noindent {\bf Counterexample to complete reducibility}:
Consider the module $V$ of dimension 3 spanned over $\k$ by
$v_{-1},v_0,v_1$ and with defining relations: 
$Kv_i = q^i v_i$; $v_0$ is annihilated by $E,X,Y,F$; $v_1$ is
killed by $E,X$; $F,Y$ kill $v_{-1}$; and, finally 
$$Fv_1 = v_{-1}, \qquad Ev_{-1} = v_1, \qquad Yv_1 = v_0,
 \qquad Xv_{-1} = -q^{-1}v_0.$$

In order to satisfy relation \eqref{r5}, we set $C_0=0$. The space $V$ and
$V_0 = \k v_0 = V(1)$ are easily seen to be $A$-modules. However, any
complement of $V_0$ in $V$ must contain a vector of the form $v = v_1 +
cv_{-1}$, and then $Yv = v_0 \in V_0$. Thus $V$ does not contain a
submodule complementary to $V_0$.

\smallskip
We also remark that the trivial module $V_0 = V(1)$ has no resolution by
Verma modules. (Such resolutions have been useful in the theory of
semisimple Lie algebras.)  For, if we had $Z(r_2) \to Z(r_1) \to V(1)$,
then $r_1 = 1$, and then $W(1) = \k[Y,F](Y,F)v_1$ would be the radical of
$Z(r_1) = Z(1)$. But then we must have $Z(r_2) \to W(1) \to 0$, whence
$r_2 = q^{-1}$ (by looking at the highest weight in both modules) and
$v_{r_2} \mapsto Yv_1$. But then the image of the map is $\k[Y,F] Yv_1$,
and $Fv_1$ is not in the image of this map.    \hfill$\Diamond$

\medskip

Recall that we write $\eps$ for $\pm 1$.
Every $K$-semisimple finite dimensional simple module is of the form
$V(r)$ for some $r = \eps q^n$, since $V_C(r)$, and hence $V(r)$, is
infinite dimensional, if $r$ is not of this form. Since char $\k \neq 2$,
every finite dimensional module is $K$-semisimple (cf. \cite[\S
2.3]{Ja}).

The main theorem of this section is the following.

\begin{theorem}\label{fd}
The simple module $V(r)$ is finite dimensional if and only if $r = \pm
q^n$ and there is a (least) integer $m>1$ so that $\alpha_{r,n-m+2} = 0$.
Furthermore, in this case,
$$V(r) = \displaystyle \bigoplus_{i=0}^{n-m} V_C(q^{-i}r)$$
\end{theorem}

\begin{proof}
Suppose $V = V(r)$ is finite dimensional simple, so $r=\eps q^n, n\in
\zz,$ as was observed above. It must be standard cyclic, so we can apply
Theorem \ref{s5} above to $V$. By \cite[Theorem 2.9]{Ja}, $V$ is a direct
sum of simple $V_C(t)$'s, each of which is finite dimensional, and
completely known, by \cite[Theorem 2.6]{Ja}.

Clearly, $v_{\eps q^{-1}}$ must be zero in $V$, else $V$ would contain a
copy of the infinite dimensional \uqsl-Verma module $Z(\eps q^{-1}) =
V(\eps q^{-1})$ (which is also simple by \cite[Proposition 2.5]{Ja}). So
let $v_{\eps q^m}$ be the ``least" nonzero \uqsl-maximal vector in $V$.
Then $v_{\eps q^{m-1}} = 0$ in $V$. Using Theorem \ref{s5} again, we can
consider $v_{\eps q^{m-1}}$ to be the image, under the quotient $\pi:
Z(r)\onto V$, of a vector $\tilde{v}_{\eps q^{m-1}}\in Z(r)$, defined as
in Theorem \ref{s5}. If $\alpha_{r,n-m+2}\not= 0$, then equation
\eqref{Rn} shows that $X\tilde{v}_{\eps q^{m-1}}$ is a non-zero multiple
of $\tilde{v}_{\eps q^m}$, so $\pi(X\tilde{v}_{\eps q^{m-1}})$ is
non-zero in $V$, which is a contradiction. Therefore,
$\alpha_{r,n-m+2}=0$.

Furthermore, $v_{\eps q^{l}}\not =0$ for $l=m,\ldots,n$. Indeed, if
$v_{\eps q^l}=0$ for some $m+1\le l\le n$, then $Xv_{\eps q^{l-1}}=0$
according to equation \eqref{Rn}, but this is a contradiction because
$v_{\eps q^n}$ is (up to a scalar) the only highest weight vector in $V$
since $V$ is simple. For the same reason, Theorem \ref{s5} implies that
$n-m+2$ is the smallest positive integer $d>1$ so that $\alpha_{r,d} =
0$.

Conversely, if there exists a $m \in \zz$ so that $\alpha_{r,n-m+2} = 0$,
then assuming that $m$ is the least such integer, we can give the \uqsl
module $V = \sum_{i=m}^n \uqsl v_{\eps q^i}$ the structure of a finite
dimensional $A$-module, using the equations worked out by Theorem
\ref{s5} and \cite[Theorem 2.6]{Ja}.
\end{proof}

We remark that one can write down the Weyl Character Formula for a simple
finite dimensional $A$-module $V$, because this formula is known for
$V_C(q^{-i}r)$.

\section{Verma modules I: Maximal vectors}\label{ssubquot}

\noindent One of the basic questions about the induced modules $Z(r)$ is:
what are their maximal vectors? The main result of this section is a step
towards a full answer to this question.

\begin{theorem}\label{s9}
We consider $Z(r)$ for any $r \in \k^\times$.
\begin{enumerate}
\item If $Z(r)$ has a maximal vector of weight $t = q^{-n}r$, then it is
unique up to scalars and $\alpha_{r,n+1} = 0$.

\item We have: $\dim_\k \Hom_A(Z(r'),Z(r)) = 0$ or $1$ for all $r,r'$,
and all nonzero homomorphisms between two Verma modules are injective.
\end{enumerate}
\end{theorem}

\noindent Part (2) follows from the first part and from the fact that
$B_- = \k[Y,F]$ is an integral domain.

Thus, a necessary condition for $Z(r)$ not to be simple (for general $r
\in \k^\times$) is that $\alpha_{r,m} = 0$ for some $m \geq 0$. Moreover,
if $r\neq\pm q^n$ ($n \in \zz$), then, from the previous section, this
condition is also sufficient, i.e. the converse to part (1) holds as
well.

To prove the first part of the theorem, we imitate \cite[Lemma 4]{Kh},
and then \cite[\S 14]{Kh}. First, we show the following lemma.

\begin{lemma}\label{compVerma}
Let $r \in\k^\times$. The following equalities hold for $v_r\in Z(r)$:
\begin{align*}
[X, F^nY^m]v_r & = -F^n \sum_{j=0}^{m-1} q^j Y^j C_0 Y^{m-1-j}v_r -
q^{m+n-1}r^{-1}\brac{q^n} F^{n-1}Y^{m+1}v_r\\
[E, F^nY^m]v_r & = -F^n \sum_{j=0}^{m-2} \brac{q^{j+1}} Y^j C_0 Y^{m-2-j}
v_r + \brac{q^n}\brac{q^{1-m-n}r} F^{n-1}Y^m v_r
\end{align*}
\end{lemma}

\begin{proof}
\begin{eqnarray*}
X(F^nY^mv_r) & = & [X, F^nY^m]v_r = [X,F^n]Y^mv_r+F^nXY^mv_r \\
 & = & [X,F^n]Y^mv_r+F^n\sum_{j=0}^{m-1}q^jY^j(XY-qYX)Y^{m-1-j}v_r \\
 & = & [X,F^n]Y^mv_r-F^n\sum_{j=0}^{m-1}q^jY^jC_0Y^{m-1-j}v_r
\end{eqnarray*}


We have to expand the first term.
\begin{eqnarray*}
[X,F^n]Y^mv_r & = & \sum_{j=0}^{n-1}F^j[X,F]F^{n-1-j}Y^mv_r \\
 & = & -\sum_{j=0}^{n-1}F^jYK^{-1}F^{n-1-j}Y^mv_r \\
 & = & -\sum_{j=0}^{n-1} q^mr^{-1}q^{2n-2-2j}F^{n-1}Y^{m+1}v_r \\
 & = & -q^mr^{-1}(\sum_{j=0}^{n-1}q^{2(n-1-j)})F^{n-1}Y^{m+1}v_r \\
 & = & -q^mr^{-1}\frac{q^{2n}-1}{q^2-1}F^{n-1}Y^{m+1}v_r \\
 & = & -q^{m+n-1}r^{-1}\brac{q^n}F^{n-1}Y^{m+1}v_r
\end{eqnarray*}

This proves the first equality of the lemma. 
We now turn to the second one. We have 
$E(F^nY^mv_r) = [E, F^nY^m]v_r = [E,F^n]Y^mv_r+F^n[E,Y^m]v_r$, 
so let us compute these two terms separately.  

\begin{eqnarray*}
[E,F^n]Y^mv_r & = & \sum_{i=0}^{n-1}F^i[E,F]F^{n-1-i}Y^mv_r \\
 & = & \sum_{i=0}^{n-1}F^i\frac{K-K^{-1}}{q-q^{-1}}F^{n-1-i}Y^mv_r \\
 & = & \sum_{i=0}^{n-1} F^i \frac{q^{-2(n-1-i)-m}r -
	q^{2(n-1-i)+m}r^{-1}}{q-q^{-1}} F^{n-1-i} Y^mv_r\\
 & = & \sum_{i=0}^{n-1} \frac{q^{-2(n-1-i)-m}r -
	q^{2(n-1-i)+m}r^{-1}}{q-q^{-1}} F^{n-1}Y^mv_r \\
 & = & \sum_{i=0}^{n-1} \frac{q^{-2i-m}r - q^{2i+m}r^{-1}}{q-q^{-1}}
	F^{n-1}Y^mv_r \\
 & = & \frac{\frac{q^{-2n}-1}{q^{-2}-1}q^{-m}r -
	\frac{q^{2n}-1}{q^2-1}q^mr^{-1}}{q-q^{-1}} F^{n-1}Y^mv_r \\
 & = & \brac{q^n}\brac{q^{1-n-m}r}F^{n-1}Y^mv_r
\end{eqnarray*}

\begin{eqnarray*}
F^n[E,Y^m]v_r & = &  F^n\sum_{j=0}^{m-1}q^{-j}Y^j
	(EY-q^{-1}YE)Y^{m-1-j}v_r \\
 & = & F^n\sum_{j=0}^{m-1}q^{-j}Y^jXY^{m-1-j}v_r\\
 & = & - F^n \sum_{i=0}^{m-1} q^{-i}Y^i \sum_{j=0}^{m-2-i}
	q^jY^jC_0Y^{m-2-j-i}v_r \\
 & = & -F^n \sum_{j=0}^{m-2} q^jY^j \left( \sum_{i=0}^{m-j}
	q^{-i}Y^iC_0Y^{-i} \right) Y^{m-2-j}v_r \\
 & = & - F^n \sum_{k=0}^{m-2} \left( \sum_{i=0}^k q^{k-2i} \right) Y^k
	C_0Y^{m-2-k}v_r \\
 & = & -F^n\sum_{k=0}^{m-2}\brac{q^{k+1}}Y^k C_0Y^{m-2-k})v_r \\
 & = & -F^n\sum_{k=0}^{m-2}\brac{q^{k+1}}Y^kC_0Y^{m-2-k}v_r
\end{eqnarray*}
This completes the proof of the lemma.
\end{proof}

\noindent {\bf Convention.}
An element $v \in Z(r) = \k[Y,F]v_r \cong \k[Y,F]$ can be viewed as a
polynomial in $Y$, with coefficients in $\k[F]$. We now define the {\it
leading term} and {\it lower order terms} of $v$ to be these terms 
with respect to the $Y$-degree.

\begin{lemma}\label{s8}
Let $r\in \k^\times$.
The following relations are valid in $Z(r)$.
\begin{enumerate}
\item Any $z \in \mf{Z}(\uqsl)$ acts on $F^n Y^mv_r$ by
$$z F^n Y^mv_r = F^n(\xi(z)(q^{-m}r) Y^mv_r + l.o.t.) \in
Z(r)_{q^{-m-2n}r}$$
\item In particular, $C_0F^nY^mv_r = F^n(c_{0,q^{-m}r}Y^m + l.o.t.)v_r$.
\item If $v \in Z(r)_{q^{-m}r}$ satisfies $Xv = 0$, then, up to scalars, we
have
$$v = (q^{m-2}r^{-1})Y^mv_r - \left( \sum_{j=0}^{m-1}
q^{m-1-j}c_{0,q^{-j}r} \right) FY^{m-2}v_r + l.o.t.$$
\item Similarly, if $v \in Z(r)_{q^{-m}r}$ satisfies $Ev = 0$, then, up to scalars,\\
\hspace*{1.6ex} (a) $v = F^n v_{\eps q^{n-1}}$, where $r = \eps
q^{m-n-1}$ (for some $n > 0$),\mbox{ OR }\\
\hspace*{1.6ex} (b) $v = \brac{q^{m-2}r^{-1}}Y^mv_r - \left( \sum_{j=0}^{m-1}
\brac{q^{m-1-j}}c_{0,q^{-j}r} \right) FY^{m-2} v_r + l.o.t.$
\end{enumerate}
\end{lemma}

\begin{proof}
\hfill
\begin{enumerate}
\item We only need to show this for the case $n=0$ because $z \in
\mf{Z}(\uqsl)$. Firstly, by weight considerations, $z \in
\text{End}_\k(Z(r)_t)$ for every $t$. Now recall that $z - \xi(z) = FUE$
for some $U \in \uqsl$. From above, $E$ takes $Y^m$ into lower order
terms, whence so does $FUE$. Therefore $zY^mv_r = \xi(z)Y^mv_r + l.o.t.$,
and $\xi(z)$ acts on the vector $Y^mv_r$ by
$\xi(z)(q^{-m}r)$, as claimed.

\item This is now obvious.

\item We first claim that any vector killed by $X$ must be of the form
$Y^m + l.o.t.$ (up to scalars). For, if 
$$v=F^n(Y^{m-2n} +
a_1 FY^{m-2n-2} + \dots)v_r = F^nY^{m-2n} + l.o.t.\ ,$$ 
then, by the above lemma, we have
$$Xv = -q^{m-n-1} r^{-1} \brac{q^n} F^{n-1}Y^{m-2n+1} + l.o.t.$$
and this is nonzero if $n>0$, because $q$ is not a root of unity. Hence
such a $v$ cannot be a solution.

Next, any solution is unique up to scalars, because given two such vectors
$v_i = Y^m + l.o.t._i$ (for $i=1,2$), we have $X(v_1-v_2) = 0$. However,
$v_1-v_2 = l.o.t._1 - l.o.t._2$, and hence must be zero from above.

Finally, from Lemma \ref{compVerma},
\begin{align*}
XY^mv_r =& -\sum_{j=0}^{m-1} q^j Y^j C_0 Y^{m-1-j} v_r\\
=& -Y^{m-1}\sum_{j=0}^{m-1} q^{m-1-j} c_{0,q^{-j}r} + l.o.t. 
\ \mbox{ (by part (2))}.
\end{align*}
Similarly, $XFY^{m-2}v_r = -q^{m-2}r^{-1}Y^{m-1} + l.o.t.$.
Hence, if $Xv=0$, then in order that the two highest degree (in $Y$)
terms cancel, $v$ must be of the given form.

\item Now suppose $Ev=0$ for some $v \in Z(r)_{q^{-m}r}$. Once again, if
$$v = F^n(Y^{m-2n} + a_1 FY^{m-2n-2} + \dots)v_r,$$ then
$$Ev = \brac{q^n} \brac{q^{1-m+2n-n}r}F^{n-1}Y^{m-2n}v_r + l.o.t.$$
and if this is zero, then $n=0$, or $r = \pm q^{m-n-1}$.

If $n=0$, then a similar analysis as above reveals that
$$EY^m = -Y^{m-2} \sum_{j=0}^{m-2} \brac{q^{1+(m-2-j)}} c_{0,q^{-j}r}
+l.o.t.,$$

\noindent and
$$EFY^{m-2} = \brac{q^{2-m}r}Y^{m-2} + l.o.t.$$ 

\noindent Therefore in this case (by the same argument as above), we must
have\\ $v = \brac{q^{m-2}r^{-1}} Y^m - bFY^{m-2} + l.o.t.$, where
$$b = \sum_{j=0}^{m-2} \brac{q^{m-1-j}}c_{0,q^{-j}r} = \sum_{j=0}^{m-1}
\brac{q^{m-1-j}}c_{0,q^{-j}r}$$

\noindent since $\brac{q^0}=0$.

On the other hand, if $n>0$, then $v = F^n v'$, where
$v' = Y^{m-2n} + l.o.t. \in Z(r)_{\pm q^{n-1}}$.
But $n>0$, so the equations of Theorem \ref{s5} apply, and we can write
$v$ as a sum of vectors in various \uqsl-Verma modules. But now,
\uqsl-theory gives us that $v = F^n v_{\pm q^{n-1}}$, where $n$ satisfies
the given conditions.
\end{enumerate}
\end{proof}

\begin{proof}[Proof of Theorem \ref{s9}]
The vector $v$ is maximal if and only if
$Ev=Xv=0$. Hence, $v$ is monic and unique
up to scalars according to the previous lemma. Using the last two parts,
we can write $v$ in two different ways.

Therefore, 
$$\brac{q^{m-2}r^{-1}} \sum_{j=0}^{m-1} q^{m-1-j}c_{0,q^{-j}r} =
(q^{m-2}r^{-1}) \sum_{j=0}^{m-1} \brac{q^{m-1-j}}c_{0,q^{-j}r}.$$

\noindent Subtracting, we get
$$\sum_{j=0}^{m-1} \bigg[ (q^{m-2}r^{-1}) \brac{q^{m-1-j}} - q^{m-1-j}
\brac{q^{m-2}r^{-1}} \bigg] c_{0,q^{-j}r} = 0.$$

\noindent Finally, using equation \eqref{s1}, we get
$$\displaystyle \sum_{j=0}^{m-1} \brac{q^{1-j}r}c_{0,q^{-j}r} = 0, \mbox{
that is, } \alpha_{r,m+1}=0.$$
\end{proof}

\section{Verma modules II: Non-integer case}\label{snonint}

\begin{stand}
From now on, unless otherwise stated, we assume that $C_0 = p(C) \neq 0$,
or $p \neq 0$. \hfill$\Diamond$
\end{stand}

Suppose $r \neq \pm q^n$ for any $n \in \Z$. Then the Verma module $Z(r)$
becomes very easy to describe. We observe that the equations \eqref{Rn},
\eqref{Sn} are valid for all $n$, so the set $\{ F^j v_{q^{-i}r} : i,j
\geq 0 \}$ is a basis for $Z(r)$.

\begin{theorem}[Non-integer power case]\label{nonint}
Suppose $r \neq \pm q^n$ for any $n \in \zz$. Then
\begin{enumerate}
\item $Z(r)$ is a direct sum of \uqsl-Verma modules $Z_C(q^{-i}r)$, one
copy for each $i$.

\item The submodules of $Z(r)$ are precisely of the form $Z(t) =
\k[Y,F]v_t$, where $t = q^{-n}r$ for every root $n$ of $\alpha_{r,n+1}$.
In particular, all these submodules lie in a chain, and $Z(r)$ has finite
length.
\end{enumerate}
\end{theorem}

\begin{proof}
The first part is a consequence of Theorem \ref{s5} and of the
observation above.  Next, if $M$ is a submodule of $Z(r)$ containing a
vector of highest possible weight $t = q^{-n}r$, then we claim that $M =
Z(t) = \k[Y,F]v_t$. To start with, $v_t \in Z(r)$ is the unique maximal
vector in $Z(r)$ of weight $t$ up to scalars, by Theorem \ref{s9} above.
Hence $v_t \in M$. We now show that $M \subset k[Y,F] v_t$.

Suppose, to the contrary, that $v \in M$ is of the form 
$$v = p(Y,F) v_t +
a_1 F^{i_1}v_{qt} + \dots + a_n F^{i_n}v_r.$$ 
We may assume that
$p(Y,F) = 0$ because $v_t \in M$. We know 
(by \cite[Theorem 2.5]{Ja}),
that the \uqsl-Verma modules $Z_C(q^it)$ are simple, so $E^lv\in
\mathrm{span}\{v_{q^it}| 1\le i\le n\}$ for some $l\gg 0$ and $E^lv\not
=0$. Therefore, since all these vectors are in different $K$-eigenspaces,
we conclude that $v_{q^it}\in M$ for some $i\ge 1$. This is a
contradiction since, by assumption, $q^it$ is not a weight of $M$ if
$i\ge 1$.
\end{proof}

\begin{remark}
In the above theorem, the successive subquotients are the simple modules
$V(t)$, and all the modules described in this section are
infinite-dimensional.
\end{remark}

\section{Verma modules III: Integer case}\label{sint}

In \S \ref{snonint}, we assumed that $r\not=\pm q^n$. 
In this section, we
treat the remaining case, namely $r = \pm q^n$. In this case, it may
happen that the simple module $V(r)$ is finite dimensional -- see Section
\ref{sfinite}. 

The main result is the following.

\begin{theorem}[Integer power case]\label{int}
Suppose $r=\pm q^n$ for $n \in \zz$. Suppose $0 = n_0 < n_1 <  \dots <
n_k \leq n+1$ are the roots to $\alpha_{r,n+1}$. Denote $t_i =
q^{-n_i}r$. Then
\begin{enumerate}
\item $Z(r)$ is a direct sum of \uqsl-Verma modules $Z_C(q^{-i}r)$ for $0
\leq i < n_k$, \emph{and} the $A$-Verma module $Z(t_k)$.
\item $Z(r)$ has the following filtration
$$Z(r) = Z(t_0) \supset W(t_0) \supset Z(t_1) \supset W(t_1) \supset
\dots \supset W(t_{k-1}) \supset Z(t_k) \supset W(t_k)$$

\noindent where the successive subquotients are, respectively,
$$V(t_0),\ V((q^3t_1)^{-1}),\ V(t_1),\ V((q^3t_2)^{-1}),\ \dots,
V(t_{k-1}),\ V((q^3t_k)^{-1}),\ V(t_k)$$

\item If $Z(t_k)$ is not simple, then it has a unique maximal submodule
of the form $Z(t)$ for some $t = \pm q^{-N}$. We then know the
composition series of $Z(t_k)$ by Theorem \ref{nonint} in this case, or
if $t_k = -1$.

\item The $V(t_i)$'s are finite dimensional (for $0 \leq i < k$).
\end{enumerate}
\end{theorem}

\noindent This theorem is similar to a corresponding one in \cite{Kh}. We
will need the following lemma.

\begin{lemma} \label{lemmaint}
Suppose that $V(r)$ is finite dimensional. If $Z(t) \hookrightarrow Z(r)$
is the largest Verma submodule in $Z(r)$ (with $t = q^m$ for some $-1
\leq m<n$), and $W(r)$ denotes the unique maximal submodule of $Z(r)$ (so
that $Z(r)/W(r)\cong V(r)$), then  $W(r)/Z(t)\cong V((q^3t)^{-1})$.
\end{lemma}

\begin{proof}
From the definition of the Casimir operator $C$, it follows immediately
that $c_r = c_{(q^2 r)^{-1}}$, whence $c_{0,r}= c_{0,(q^2r)^{-1}}$. We
claim that $\alpha_{r,2n+4}=0$. Indeed,
\begin{eqnarray*}
\alpha_{r,2n+4} & = & \sum_{j=0}^{2n+2}\brac{q^{1-j+n}}c_{0,q^{-j+n}} \\
 & = & \sum_{j=-n-1}^{n+1}\brac{q^j}c_{0,q^{j-1}} \\
 & = & \sum_{j=1}^{n+1}(\brac{q^j}c_{0,q^{j-1}} +
	\brac{q^{-j}}c_{0,q^{-j-1}}), \mbox{ since } \brac{q^0} = 0\\
 & = & \sum_{j=1}^{n+1} \brac{q^j}(c_{0,q^{j-1}} - c_{0,q^{-j-1}}),
	\mbox{ since } \brac{a} + \brac{a^{-1}} = 0\ \forall a \\
 & = & \sum_{l=0}^n \brac{q^{l+1}}(c_{0,q^l} - c_{0,q^{-l-2}})\\
 & = & 0 \qquad \mbox{ from above}.
\end{eqnarray*}

Thus, $t$ is the first root after $r$ for $\alpha_{r,n}$, 
if and only if
$(q^3r)^{-1}$ is the first root after $(q^3t)^{-1}$ for
$\alpha_{(q^3t)^{-1},n}$.

But now, the quotient $W(r) / Z(t)$ has a vector of highest weight
$(q^3t)^{-1}$: if $v_{qt}$ is the lowest \uqsl-maximal vector in $V(r)$,
and $t = \eps q^m$, then $F^{m+2}v_{qt}$ is \uqsl-maximal and of highest
weight in the quotient. But it has weight $q^{-2m-4}\eps q^{m+1} =
(q^3t)^{-1}$ as claimed.

Thus, $W(r) / Z(t)$ has a subquotient of the form $V((q^3t)^{-1})$. But
one can check that they have the same characters. Hence they are equal.
\end{proof}

\begin{proof}[Proof of Theorem \ref{int}]
To simplify the notation, let us assume that $r=q^n$; the case 
$r=-q^n$ is similar. Suppose that the simple module $V(r)$ is 
finite dimensional. Then, by Theorem \ref{fd}, $\alpha_{r,n-m+2}=0$ 
for some $m>1$ and $n-(m-1)>0$. We assume that $m$ is the smallest 
such integer. Then, from the proof of that same theorem, we know 
that $v_{q^{m-1}}$ is maximal in $Z(r)$. Therefore, setting 
$n_1=n-(m-1)$ and $t_1=q^{-n_1}r=q^{m-1}$, we get that $Z(t_1)\into Z(r)$.

We can repeat the same procedure with $Z(t_1)$. If its simple top 
quotient $V(t_1)$ is finite dimensional, then there exists a 
(smallest) integer $m_1$ such that $\alpha_{(m-1)-m_1+2}=0$ 
for some $m_1>1$ and $(m-1)-(m_1-1)>0$. Again, $v_{q^{m_1-1}}$ 
is maximal in $Z(t_1)$, so $Z(t_2)\into Z(t_1)$ where 
$t_2=q^{-n_2}t_1=q^{m_1-1}$ and $n_2=m-1-(m_1-1)$. Note 
that $n-(m_1-1)=n-(m-1)+(m-1)-(m_1-1)>0$.

We can continue repeating this procedure and get a chain of 
Verma submodule $Z(r)\supset Z(t_1)\supset Z(t_2)\supset\ldots$.
Set $n_0=n,m_0=m,d_0=n-(m-1)+1$ and $d_i=(m_{i-1}-1)-m_i+2$ for 
$i\ge 1$. Since $n-(m_i-1)>0$ for all $i\ge 1$ (as noted in the 
previous paragraph for $i=1$), this procedure must stop for 
some positive integer $k$. This means that $Z(t_k)\subset Z(r)$, 
but the top quotient of $Z(t_k)$ is not finite dimensional.

Using Theorem \ref{s5}, this proves part (1). We can now apply Lemma
\ref{lemmaint} to each successive inclusion $Z(t_i)\subset Z(t_{i-1})$,
and part (2) is proved. Part (4) follows from the first two parts.

It remains to show part (3); namely, that $W(t_k) = 0$ or $Z(t)$ for some
$t$. So suppose $Z(t_k)$ is not simple. Let $v_t$ be the highest possible
maximal vector in $Z(t_k)$, that is not of weight $t_k$ (i.e. it has
``smaller" weight). Thus $t = q^{-n}t_k$ for some $n$, and $v_t = Y^n
v_{t_k} + l.o.t.$, from Lemma \ref{s8} above.

Now, any weight vector $v \in W(t_k)$ is of the form $g(Y,F)v_t + F^l
h(Y,F)v_{t_k}$, where $h$ is {\it monic}. (This follows from the
Euclidean algorithm for polynomials $(\k[F])[Y]$, because $v_t$ is
monic.) Further, $l>0$, since we are not considering the case $t_k = \pm
q^{-1}$, which we know by Section \ref{snonint}.

To show $W(t_k) = Z(t)$, we must prove that $h = 0$ for each such $v_x$.
Suppose not. Let $v_x \in W(t_k)$ be a weight vector of highest possible
weight $x$, such that $h \neq 0$. Now, $Ev_x \in W(t_k)$, so by
maximality of $x$, $Ev_x \in Z(t) = \k[Y,F]v_t$, hence $E(v_x -
g(Y,F)v_t) \in Z(t)$. Hence, we get that $E(F^l h(Y,F)v_{t_k}) \in Z(t)$.

This is in the \uqsl-span of $v_{t_k}, Yv_{t_k}, Y^2v_{t_k}, \dots,
Y^{n-1}v_{t_k}$, so if it is in $Z(t)$, then it must be zero, by the PBW
Theorem. Hence $E(F^l h(Y,F) v_{t_k}) = 0$. But now, part (4) of
Lemma \ref{s8} above, gives us that $F^l h(Y,F) v_{t_k} = F^l
v_{\pm q^{l-1}}$.

Hence we finally get that $v' = F^l v_{\pm q^{l-1}} \in W(t_k)$. Hence
$X^l v' \in W(t_k)$. From the following lemma, this means that (up to a
nonzero scalar,) $v_{\pm q^{-1}} \in W(t_k)$. But $t$ was ``lower" than
$\pm q^{-1}$ from above, hence this is a contradiction, and no such $v_x$
exists.
\end{proof}

\begin{lemma}\hfill
\begin{enumerate}
\item $[F^{j+1},X] = q^j \brac{q^{j+1}}F^jYK^{-1}$
\item If $r = \eps q^n$, then $F^{j+1}v_{\eps q^j}$ is \uqsl-maximal (for
each $-1 \leq j \leq n$), and\\
$X(F^{j+1}v_{\eps q^j}) = -\brac{\eps q^{j+1}} F^jv_{\eps q^{j-1}}$.
\end{enumerate}
\end{lemma}

\begin{proof}
For the first part, we compute, using the defining relations:
\begin{align*}
[F^{j+1}, X] & = \sum_{i=0}^j F^{j-i} [F, X] F^i = \sum_{i=0}^j F^{j-i} Y
K^{-1} F^i = \sum_{i=0}^j F^{j-i} Y q^{2i} F^i K^{-1}\\
& = \sum_{i=0}^j q^{2i} \cdot F^j Y K^{-1} = \frac{q^{2j+2} - 1}{q^2 - 1}
F^j Y K^{-1} = q^j \brac{q^{j+1}} F^j Y K^{-1}
\end{align*}

\noindent as claimed. Next, suppose $r = \eps q^n$ for some $n$. We then
compute:
\begin{align*}
E \cdot F^{j+1} v_{\eps q^j} & = F^{j+1} \cdot E v_{\eps q^j} +
\sum_{i=0}^j F^{j-i} [E,F] F^i v_{\eps q^j}\\
& = 0 + \sum_{i=0}^j F^{j-i} \frac{K - K^{-1}}{q - q^{-1}} F^i v_{\eps
q^j}\\
& = \sum_{i=0}^j F^{j-i} \frac{\eps}{q - q^{-1}} (q^{j-2i} - q^{2i-j})
F^i v_{\eps q^j}\\
& = \frac{\alpha \eps}{q - q^{-1}} F^j v_{\eps q^j},
\end{align*}

\noindent where $\alpha = \sum_{i=0}^j q^{j - 2i} - q^{2i - j} = (q^j +
q^{j-2} + \ldots + q^{-j}) - (q^{-j} + q^{-j+2} + \ldots + q^j) = 0$.
Thus $F^{j+1} v_{\eps q^j}$ is $\uqsl$-maximal as claimed.

Finally, we show the last assertion, for which we need the first part of
this lemma, as well as equations \eqref{Rn}, \eqref{Sn}. We compute:
\begin{align*}
XF^{j+1} v_{\eps q^j} & = F^{j+1} X v_{\eps q^j} - q^j \brac{q^{j+1}} F^j
Y K^{-1} v_{\eps q^j}\\
& = F^{j+1} \cdot \left( -\frac{\alpha_{r, n-j+1}}{\brac{\eps q^{j+2}}}
v_{\eps q^{j+1}} \right) - q^j \brac{q^{j+1}} \eps q^{-j} F^j Y v_{\eps
q^j}\\
& = -\frac{\alpha_{r, n-j+1}}{\brac{\eps q^{j+2}}} F^{j+1} v_{\eps
q^{j+1}} - \brac{\eps q^{j+1}} F^j Y v_{\eps q^j}\\
& = -F^j \brac{\eps q^{j+1}} \left( Y v_{\eps q^j} + \frac{\alpha_{r,
n-j+1}}{\brac{\eps q^{j+1}} \brac{\eps q^{j+2}}} F v_{\eps q^{j+1}}
\right)\\
& = -\brac{\eps q^{j+1}} F^j (Y v_{\eps q^j} + d_{r, n-j+1} F v_{\eps
q^{j+1}})\\
& = -\brac{\eps q^{j+1}} F^j v_{\eps q^{j-1}}
\end{align*}

\noindent and we are done.
\end{proof}

\section{Category \calo}\label{sfinl}

Our goal in this section is to show that the category \calo (defined
below) is a highest weight category in the sense of \cite{CPS} and that
it can be decomposed into a direct sum of subcategories (``blocks''),
each of which contains only finitely many simple modules
\footnote{The original paper \cite{BGG1} achieves this using the
eigenvalues of the Casimir operator. However, unlike their case, we will
see later that our algebra $A$ has trivial center (if $C_0\neq 0$).
Therefore, such an approach fails and we have to do more work (similarly
to \cite{Kh}).}.
We retain our assumption that $C_0 \neq 0$.

\begin{definition}
The category \calo consists of all finitely generated
$A$-modules with the following properties:
\begin{enumerate}
\item The $K$-action is diagonalizable with finite dimensional weight
spaces.
\item The $B_+$-action is locally finite.
\end{enumerate}
\end{definition}

Given $r \in \k^\times$, we claim there exist only finitely many $t =
q^{-n}r$ such that $Z(t) \hookrightarrow Z(r)$. If we have such an
embedding, then $\alpha_{r,n+1}=0$ by Theorem \ref{s9}, so we have to see
that this is true for only finitely many $n$ if $r$ is fixed. Proposition
\ref{polynom1} says that $\alpha_{r,n+1}$ is a nonvanishing function,
multiplied by a polynomial in $q^{n+1}$, if $r$ is fixed. This polynomial
can be factored as $\prod_{i=1}^{L}(q^{n+1}-z_i)$ where $z_1, \ldots,
z_L$ are the roots of the polynomial. This will be zero only for values
of $n$ such that $q^{n+1}=z_i$ for some $i$; since $q$ is assumed not to
be a root-of-unity, there are only finitely many such $n$.

We claim also that, fixing $r$, there are only finitely many $s$ of the
form $s=q^nr$ such that $Z(r)\into Z(s)$. This is because, if we have
such an embedding, then $\alpha_{s,n+1}=0$ by Theorem \ref{s9} and since
for fixed $r,\ \alpha_{s,n+1} = \alpha_{q^n r,n+1}$ is (as above)
essentially a polynomial in $q^{n+1}$ (by looking at the expansion of
$b(S,T)$ in Proposition \ref{polynom1}), it vanishes for only finitely
many values of $n$.

Let us fix $r$ and consider the {\it maximal} $n \geq 0$ so that
$\alpha_{q^nr,n+1} = 0$. That such an $N$ exists follows from the
observation (in the previous paragraph) that the set of such $n$ is
finite.  Set $r_0 = q^Nr$, so $\alpha_{r_0,N+1}=0$.

Define $S(r)$ to be the set of all $t = q^{-m}r_0$, so that
$\alpha_{r_0,m+1} = 0$. This is a finite set.

We now introduce a graph structure on $\k^\times$ by connecting $t$ and
$r$ by an edge if and only if
$Z(r)$ has a simple subquotient $V(t)$ or $Z(t)$ has a
simple subquotient $V(r)$. The component of this graph containing $r$ is
denoted $T(r)$.

\begin{proposition}\label{finite}
\hfill
\begin{enumerate}
\item If $t\in S(r)$, then $S(t)=S(r)$.

\item For each $r\in\k^{\times}$, $T(r) \subset S(r)$. In particular,
$T(r)$ is finite for each $r$.

\item Every Verma module has finite length.
\end{enumerate}
\end{proposition}

\begin{proof}
The proof of part (1) is in two parts. The first one is the following
equality:
\begin{equation}\label{eqn}
\alpha_{q^nr,n+m+1}=\alpha_{q^nr,n+1}+\alpha_{r,m+1}
\end{equation}
We provide a proof of this equality using the definition of $\alpha$:
\begin{eqnarray*}
\alpha_{q^nr,n+m+1} & = &
	\sum_{j=0}^{n+m-1}\brac{q^{1-j}(q^nr)}c_{0,q^{-j}(q^nr)} \\
 & = & \sum_{j=0}^{n-1} \brac{q^{1-j} (q^nr)}c_{0,q^{-j}(q^nr)} +
	\sum_{j=n}^{n+m-1} \brac{q^{1-j}(q^nr)} c_{0,q^{-j}(q^nr)} \\
 & = & \alpha_{q^nr,n+1} + \sum_{i=0}^{m-1}\brac{q^{1-i-n}(q^nr)}
	c_{0,q^{-i-n}(q^nr)} \\
 & = & \alpha_{q^nr,n+1}+\sum_{i=0}^{m-1}\brac{q^{1-i}r}c_{0,q^{-i}r}\\
 & = & \alpha_{q^nr,n+1}+\alpha_{r,m+1}
\end{eqnarray*}


Now suppose $t\in S(r)$, so $t=q^{-l}r_0$ and $\alpha_{r_0,l+1}=0$. We
define $t_0$ similarly to $r_0$, so, in particular, $t_0=q^T t$ and
$\alpha_{t_0,T+1}=0$. We claim that $t_0=r_0$, which implies that
$S(t)=S(r)$. Note that, by the maximality of $t_0$, $l\le T$. We have to
show that $l=T$, so assume that, on the contrary, $l<T$. Equation
\eqref{eqn} along with $\alpha_{t_0,T+1}=\alpha_{r_0,l+1}=0$ implies that
$\alpha_{t_0,T-l+1}=0$. This last equality, now in conjunction with
$\alpha_{r_0,N+1}=0$ and equation \eqref{eqn}, yields
$\alpha_{t_0,T-l+N+1}=0$. Since $t_0=q^{N+T-l}r$ and $N+T-l>N$, this
contradicts the maximality of $N$. Therefore, $T=l$.

The proof of part (2) is also in two steps. First, we need the following
observation: Theorems \ref{nonint} and \ref{int} state that if $V(t)$
is a subquotient of $Z(r)$, then $t=q^{-m}r$ for some root $m$ of
$\alpha_{r,m+1}$. The second step consists in showing that $S(t)=S(r)$ if
$V(t)$ is a subquotient of $Z(r)$. From part (1), it is enough to show
that $t \in S(r)$. If $V(t)$ is a simple subquotient of $Z(r)$ with $t =
q^{-m}r$, then $\alpha_{r,m+1}=0$, and combining this with
$\alpha_{r_0,N+1} = 0$ and equation \eqref{eqn}, we get
$\alpha_{r_0,m+N+1} = 0$. Since $t = q^{-m-N} r_0$, this means exactly
that $t \in S(r)$.

The general case of an arbitrary $t \in T(r)$ follows from the specific
case that we just considered. 

Part (3) is a consequence of part (2) and of the fact that the simple
quotients of a Verma module occur with finite multiplicities, which, in
turn, is a consequence of the fact that the weight spaces of every Verma
module are finite dimensional.
\end{proof}

\begin{definition}\label{sf}
A finite filtration $M=F^0\supset F^1\supset\ldots F^r=\{0\}$ of a module
$M \in \calo$ is said to be \emph{standard} if $F^i/F^{i+1}$ is a Verma
module for all $i$.
\end{definition}

We construct some useful modules which admit such a filtration.
Let $a\in\k^\times, l \in \zz$;
define $Q(l)$ to be the $A$-module induced from the $B_+$-module
$B_+/N_+^l$,
and define $Q(a,l)$ to be the $A$-module induced from the $B_+$-module
$\mathcal{B}_{a,l} := B_+/((K-a), N_+^l)$,
so $Q(a,l)=A\otimes_{B_+}\mathcal{B}_{a,l}$ is a quotient of
$Q(l)=A\otimes_{B_+}B_+/N_+^l$.
Notice that the modules $Q(a,l)$ all have standard filtrations, because
$N_+^j\mathcal{B}_{a,l}/N_+^{j+1}\mathcal{B}_{a,l}$ is a $B_+$-module
on which $N_+$ acts trivially, and $\k[K,K^{-1}]$ semisimply.

Moreover, given any module $M \in \calo$ and a weight vector $m \in M$ of
weight $a$, there exists a nonzero homomorphism $f : Q(a,l) \to M$ for
some $l$, taking $\bar{1}$ to $m$, where $\bar{1}$ is the generator
$1\otimes 1$ in $Q(a,m)$. This is because $N_+$ acts nilpotently on $m$.

\begin{proposition}\label{quotient}
Every module in \calo is a quotient of a module which admits a standard
filtration.
\end{proposition}

\begin{proof}
Let $M$ be an arbitrary module in \calo. Since $M$ is finitely generated
over $A$, $M$ is Noetherian according to Corollary \ref{noeth}. Choose a
nonzero weight vector $m_1$ in the weight space $M_{a_1}$ (for some $a_1
\in \k^\times$), and an arbitrary non-zero homomorphism
$f_1:Q(a_1,l_1)\to M$ for some $l_1$, and set $N_1=\im(f_1)$. If $N_1\neq
M$, choose another homomorphism $f_2:Q(a_2,l_2)\to M$ such that
$N_1\subsetneq N_1+N_2$, where $N_2 = \im(f_2)$ (this is possible by the
remark above).

Repeating this procedure, we get an increasing chain of submodules
$N_1\subsetneq N_1+N_2\subsetneq N_1+N_2+N_3\subsetneq\ldots$ which must
stabilize since $M$ is Noetherian. This implies that
$M=N_1+N_2+\ldots+N_k$ for some $k$. It was observed above that
$Q(a_1,l_1)\oplus\cdots\oplus Q(a_k,l_k)$ has a standard filtration.
\end{proof}

\begin{proposition}\label{fl}
Every module in \calo has finite length.
\end{proposition}

\begin{proof}
This is an immediate consequence of Proposition \ref{quotient} and of the
fact that Verma modules have finite length - see Proposition
\ref{finite}.
\end{proof}

We introduce the following partial order on $\k$: $t \le s$
if and only if $t=q^l s$ for some $l \in \mathbb{Z}_{\le 0}$.

\begin{proposition}\label{extvan}
\hfill
\begin{enumerate}
\item If $s\not\in T(r)$, then $\Ext(V(r),V(s))=0$ and
$\Ext(Z(r),Z(s))=0$.
\item Simple modules have no self-extensions.
\end{enumerate} 
\end{proposition}

We omit the proof of the preceding proposition, which is same as the
corresponding statements for $H_f$; see \cite[Theorem 4]{Kh}. We only
need to show the existence of a ``good" duality functor $\F$ as in
\cite[\S 2]{Kh}. We do this now.


\remark\label{antiinv} It is easy to check that the following define an 
anti-involution $\inv$ of $A$:
$$\inv(E)=-FK,\quad \inv(F)=-K^{-1}E,\quad \inv(K)=K,\quad
\inv(K^{-1})=K^{-1},$$
$$\quad \inv(X)=Y, \quad \inv(Y)=X.$$

\noindent Note, in particular, that $\inv(C)=C$.
\medskip


\begin{definition}
Define the {\it duality functor} $\F$ on $\calo$ as follows:
if $M \in \calo$, then let $\F(M)$ be
the linear span of all $K$-semisimple vectors in $M^*$. The $A$-module
structure on $\F(M)$
is given, using the anti-involution $\inv$, by
\[ (a m^*)(m) := m^*(\inv(a)m) \qquad \mbox{ for all }
a \in A,\ m \in M,\ m^* \in
\F(M). \]
\end{definition}

As in \cite[Proposition 2]{Kh}, the duality functor $\F$ is exact,
contravariant, takes simple objects in $\calo$ to themselves, and
preserves the formal characters and the set of composition factors, of
any (finite length) object in $\calo$.
\medskip

\begin{definition}
Define $\calo(r)$ to be the subcategory of all the modules whose simple
subquotients $V(t)$ satisfy $t\in T(r)$.
\end{definition}

\begin{corollary} \label{blockd}
We have a decomposition $\calo=\bigoplus_{r} \calo(r)$.
\end{corollary}

\begin{proof}
This is an immediate consequence of the vanishing of the $\Ext$ in the
previous proposition.
\end{proof}

\begin{proposition}\label{proj}
The category \calo has enough projective objects.
\end{proposition}

\begin{proof}
Consider a component $T(r)$: we know it is finite. Pick
$s\in\k^{\times}$. Since $T(r)$ is finite, there exists an integer $n_s$
such that $N_+^{n_s}v=0$ for any $v\in M_s$ and any $M\in \calo(r)$,
where $M_s$ is the weight space of $M$ of weight $s$. For each such $s$,
choose such an $n_s$.

Since $\Ext(Z(r_1),Z(r_2))=0$ if $T(r_1)\neq T(r_2)$ according to
Proposition \ref{extvan}, it follows that $Q(s,n_s)$ decomposes as a
direct sum $Q(s,n_s)=\oplus_{r} Q(s,n_s,r)$, where $Q(s,n_s,r)$ is a
submodule, all of whose successive subquotients are in $T(r)$. It should
be noted that $Q(s,n_s,r)=0$ if $s$ is not of the form $s=q^lr$. Set
$P(s,r)=Q(s,n_s,r)$.

We claim that $P(s,r)$ is projective in \calo. Indeed,
$\Hom_\calo(P(s,r),V(t))=0$ if $t\not\in T(r)$, and if $M \in \calo(r)$,
then $\Hom_\calo(P(s,r),M)= \Hom_\calo(Q(s,n_s),M) = M_s$. Since the
$K$-action on $M$ is diagonalizable, $M\rightarrow M_s$ is an exact
functor from \calo to the category of vector spaces. Therefore, $P(s,r)$
is projective.
\end{proof}

Let $V(s)$ be a simple $A$-module. Then $P(s,s)$ admits an epimorphism
onto $V(s)$. Since $P(s,s)$ has finite length, we can express it as the
direct sum of finitely many indecomposable projective modules. This
implies that there exists an indecomposable direct summand $P(s)$ with a
non-zero homomorphism $P(s)\to V(s)$. This module $P(s)$ also admits a
standard filtration since it is a direct summand of a module with such a
filtration \cite{BGG1}.

\begin{proposition}\label{hwc}
The category \calo is a highest weight category.
\end{proposition}

\begin{proof}
The only two points that we have to prove are the following:
\begin{enumerate}
\item If $V(t)$ is a subquotient of $Z(r)$, then $t\le r$ and the
multiplicity $[Z(r) : V(r)]$ of $V(r)$ as a subquotient of $Z(r)$ is one.
\item If $Z(r)$ appears as a subquotient in a standard filtration of
$P(s)$, then $s\le r$. Moreover, $Z(s)$ appears exactly once in any such
filtration.
\end{enumerate}

The statement (1) is a consequence of the observation that if $t$ is a
weight of $Z(r)$, then $t\le r$. Moreover, the weight space of $r$ in
$Z(r)$ has dimension one.

The second part follows from the following vanishing result: if $r\not
<s$, then we have $\Ext(Z(r),Z(s))=0$, which can be proved exactly as the
analogous result (Proposition 4) in \cite{Kh}. Another approach (see
\cite[Proposition 11]{Kh}) is to use  the construction of $P(s)$ as a
direct summand of $P(s,s)$.
\end{proof}

\section{Block decompositions in highest weight
categories}\label{sblockdec}

In a highest weight category like \calo, it is possible to define a block
decomposition in several different ways. We now show why all these ways
yield the same decomposition. To begin with, we can define a
decomposition using Verma modules; this is exactly the one given by the
sets $T(r)$, and we rephrase this as the condition
\[ G_{VZ}: ``V(t) \mbox{ is a subquotient of } Z(r)."\]

\noindent Thus, $S_{VZ}(r)$ is the graph component of $\k^\times$
containing $r$, where we join $r$ and $t$ by an edge if $V(t)$ is a
subquotient of $Z(r)$, or $V(r)$ is a subquotient of $Z(t)$.

Recall that there exists an exact contravariant duality functor $\F$ (as
mentioned above) that takes simple objects to themselves, and preserves
the set of composition factors of any object of finite length in \calo.

Using this functor $\F$, we can now define $A(r) = \F(Z(r))$ and $I(r) =
\F(P(r))$ to be the co-standard and (indecomposable) injective modules
respectively. In a highest weight category like \calo, every projective
module has a standard filtration as above, and BGG Reciprocity also holds
(cf. \cite{CPS2}). In other words, $[P(r) : Z(t)] = [Z(t) : V(r)]$ for
all $t,r$.

\noindent {\bf Definitions}:
\begin{enumerate}
\item We define the property $G_{ab}$ by
\[ G_{ab} :\ ``a(t) \mbox{ is a subquotient of } b(r)" \]

\item Given $a,b$ as above, we introduce a graph structure on $\k^\times$
as follows: connect $r$ and $t$ by an edge $r$---$t$ if $G_{ab}$ holds
for the pair $(r,t)$ or $(t,r)$. Under this structure, we define the
connected component of $\k^\times$ containing $r$, to be the block
$S_{ab}(r)$.

\item We also have the {\it categorical definition of linking}: We say
$r$ and $t$ are {\it linked} if there is a chain of {\it indecomposable}
objects $V_i \in \calo$ and nonzero maps $f_i \in \hhom(V_{i-1},V_i)$
such that
$$V_0 = V(t) \mapdef{f_1} V_1 \mapdef{f_2} \dots \mapdef{f_n} V_n =
V(r)$$

\item We now define the final graph structure on $\k^\times$ as follows:
$B(r)$ is the connected component of $\k^\times$ containing $r$, where
edges denote linked objects.
\end{enumerate}

\noindent We remark that the $V_i$'s need to be indecomposable, otherwise
any two objects of \calo are linked by $0 \to M \to M \oplus N \to N \to
0$. Also note that the definition of linking is clearly symmetric, using
the duality functor $\F$.

We now explain why certain block decompositions of $\k^\times$ are the
same. Using the duality functor, it is easy to see that the conditions
$G_{VZ}$ and $G_{VA}$ are the same; hence we have $S_{VZ}(r) = S_{VA}(r)$
for all $r$. Similarly, $G_{VP}$ and $G_{VI}$ are equivalent, as are
$G_{ZP}$ and $G_{AI}$. We now have the following result.

\begin{theorem}\hfill
\begin{enumerate}
\item $T(r) = S_{VZ}(r) \subset S_{ZP}(r) \subset S_{VP}(r) \subset B(r)$
for all $r$.
\item In fact, $T(r) = B(r)$.
\end{enumerate}
\end{theorem}

\noindent In particular, all the various block classifications coincide
in our case.

\begin{proof}\hfill
\begin{enumerate}
\item We proceed step by step. The philosophy is to show that two
vertices connected by an edge in a graph structure are connected in a
bigger one.

Suppose $r,t$ are edge-connected in $T(r) = S_{VZ}(r)$. Then, by BGG
reciprocity, $[P(r): Z(t)] = [Z(t) : V(r)] >0$. Hence $S_{VZ} \subset
S_{ZP}$.
Next, if $P(r)$ has a subquotient $Z(t)$, then it clearly has a
subquotient $V(t)$ as well. Thus $S_{ZP} \subset S_{VP}$.

Finally, suppose $P(r)$ has a simple subquotient $V(t)$. We show that
$r$ and $t$ are linked. We have a sequence of maps
\[ 0 \to N \to M \to V(t) \to 0, \qquad N \hookrightarrow M
\hookrightarrow P(r) \twoheadrightarrow V(r) \]

\noindent Since $V(t)$ is indecomposable, we can choose $M$ to be
indecomposable as well.
Since \calo is finite length, we have $V(s) \hookrightarrow N
\hookrightarrow M$ for some $s$. Hence, using the duality functor $\F$,
we now have the following sequence of maps linking $V(t)$ and $V(r)$:
$$V(t) \cong \F(V(t)) \hookrightarrow \F(M) \twoheadrightarrow \F(V(s))
\cong V(s) \hookrightarrow M \hookrightarrow P(r) \twoheadrightarrow
V(r)$$

\noindent and we are done.

\item Since $T(r)$ is finite for all $r$, results from the previous
section tell us that $\calo = \bigoplus_r \calo(r)$, where each
$\calo(r)$ is a highest weight category.
Now, suppose $V(r)$ and $V(t)$ are linked via a chain $V_0 = V(r) \to V_1
\to \dots \to V_n = V(t)$. Since each $V_i$ is indecomposable, it is in a
unique block $\calo(s)$. However, since there are nonzero homomorphisms
in between successive $V_i$'s, all the $V_i$'s are in the same block. In
particular, $V(t) \in \calo(r)$, so $B(r) \subset T(r)$ for each $r$, and
hence is finite as well. Combining this with part (1) yields $T(r) =
B(r)$.
\end{enumerate}
\end{proof}


\section{Semisimplicity}\label{sweyl}

As we saw in \S \ref{sfinite}, finite dimensional $K$-semisimple
$A$-modules are not $A$-semisimple for some $A$ (or $C_0 = p(C)$).
However, we have the following result, that tells us when the result
{\it does} hold.

\begin{theorem}\label{weyl}
The following are equivalent:
\begin{enumerate}
\item For each $n \in \zz$ (and each $r = \eps q^n),\ T(r) \subset \{
r,s, (q^3s)^{-1}, (q^3r)^{-1} \}$, where $s = \eps q^m$ for some $m + 1
\in \zz$.
\item For each $n \in \zz$ (and each $r = \eps q^n$), the equation
$\alpha_{r,m} = 0$ has at most one root $m$ satisfying $2 \leq m \leq
n+1$.
\item For each $r\in\k^{\times}$, there is at most one finite dimensional
simple module in $\calo(r)$ (up to isomorphism).
\item Every finite dimensional $A$-module is completely reducible.
\end{enumerate}
\end{theorem}

\begin{proof}
Since char $\k \neq 2$, every finite dimensional module is
$K$-semisimple, and hence it is an object of the category $\calo =
\bigoplus_r \calo(r)$.

Note from Lemma \ref{lemmaint} that $t \in S(r)$ if and only if
$(q^3t)^{-1} \in
S(r)$, for every $r$ of the form $\eps q^n,\ n\in \zz$. This explains the
structure of the set in (1) above (since $T(r) \subset S(r)$).

Next, observe that the set of isomorphism classes of simple modules in a
block $T(r)$ is in bijection with the block $T(r)$ itself (under the map
$t \mapsto V(t)$ for each $t \in T(r)$), and every simple module $V(r)$
is actually in a block $T(r)$. We first show that the first three
assertions are equivalent.

If (1) holds, then any simple module in $\calo(r)$ is one of the
following:
$$V(r),\ V(s), V((q^3r)^{-1}),\ V((q^3s)^{-1}).$$

\noindent Since only $V(r)$ is finite dimensional in the above list, and
that too only when $r = \eps q^n$ for $n \in \zz$, (1) implies (3).

Similarly, if (1) does not hold then $S(r)$ contains $t_i = \eps q^{m_i}\
(i=0,1,2)$, where we assume without loss of generality
that $m_0 = n > m_1 > m_2$. Then
$\alpha_{r,m}$ has two roots by definition of $S(r)$, so (2) does not
hold either. In other words, (2) implies (1).

Now suppose that (2) does not hold. Thus there are at least two roots of
$\alpha_{rm}$. By Theorem \ref{s5}, there are weight vectors $v_{t_1},
v_{t_2}$, say of weights $t_i = \eps q^{m_i}$, in the Verma module
$Z(r)$, with $-1 \leq m_2 < m_1 \leq m_0 = n$. But then by Theorem
\ref{fd}, there are (at least) two nonisomorphic finite dimensional
simple modules, namely $V(r) = V(t_0)$ and $V(t_1)$ in $\calo(r)$, by
Theorem \ref{int} and equation \eqref{eqn}. Thus (3) does not hold
either, meaning that (3) implies (2), and the first three assumptions are
shown to be equivalent.

Now suppose (3) holds. We show complete reducibility. Suppose $M$ is a
finite dimensional $A$-module. Since the category \calo splits up into
blocks, we have $M = \bigoplus_r M(r)$, where $M(r) \in \calo(r)$. Each
$M(r)$ is finite dimensional, hence so are all its subquotients. Hence
by assumption, all subquotients of $M(r)$ are of the form $V(r)$. Since
$V(r)$ has no self-extensions in \calo by Proposition \ref{extvan}, this
shows (using \cite[Proposition A.1]{Kh} and induction on length, for
instance) that $M(r)$ is actually a direct sum of copies of $V(r)$. Hence
$M(r)$ is semisimple.

Finally, suppose (1) does not hold; we show that (4) does not hold
either. As in Theorem \ref{int}, let $r = \eps q^n$, and let $0 = n_0 <
n_1 < \dots < n_k \leq n+1$ be the various roots of $\alpha_{r,m+1}$.
Since (1) fails, we have $k \geq 2$.

Given $i \geq j$, we now define the module $W(i,j)$ to be the
$A$-submodule generated by $\{ F^{b+1}v_{\eps q^b} : n - n_i \leq b \leq
n - n_j \}$ and $Z(t_i)$. For example, $W(i,i) = Z(t_i)$ is a Verma
module, and $W(i+1,i) = W(t_i)$ is its unique maximal submodule.

We now consider the filtration
$$Z(t_0) = W(0,0) \supset W(t_0) = W(1,0) \supset W(2,0)$$

\noindent This gives a short exact sequence
$$0 \to W(1,0) / W(2,0) \to Z(t_0) / W(2,0) \to Z(t_0) / W(1,0) \to 0$$

\noindent or, in other words,
$$0 \to V(t_1) \to Z(t_0) / W(2,0) \overset{\varphi}{\longrightarrow}
V(t_0) \to 0$$

\noindent The middle term is thus a finite dimensional module of length
2. We claim that there does not exist a splitting of the map $\varphi$.
This is easy to show: any complement to $V(t_1)$, if it exists, is also
$K$-semisimple, and hence contains the highest weight vector $v_{t_0}$.
But $v_{t_0}$ generates the entire module $Z(t_0) / W(2,0)$, so there
cannot exist a complement, and (4) fails, as claimed.
\end{proof}

\begin{remark}
Note that the condition (2) above depends on the polynomial $p$, or in
other words, on the central element $C_0 = p(C)$, by means of the
polynomial $\alpha_{r,m}$. Furthermore, there {\it are} central elements
$C_0$ in \uqsl, that satisfy the condition above. We give such an example
now. (Note that in the case $C_0 = 0$, complete reducibility was
violated.) \hfill$\Diamond$
\end{remark}\hfill

\noindent {\bf Example of complete reducibility:}

\begin{stand}
For this example, $q$ is assumed to be transcendental over $\mathbb{Q}$.
\end{stand}

Take $p(C) = C_0 = (q-q^{-1})^3 C - (q-q^{-1})\qdiff$. (Note that if $A$
satisfies complete reducibility (as above) for this $p$, then it does so
for any scalar multiple of $p$.) We now show that the only finite
dimensional simple module is $V(q) = V_C(q)$, of dimension 2 over $\k$.

Let us first calculate $\alpha_{r,m}$, using the (computations in the)
proof of Proposition \ref{polynom1}. Clearly, we have
$h(qT) = [(qT + q^{-1}T^{-1}) - \qdiff](q - q^{-1})$, so \[ g(qT) = h(qT)
\brac{qT} = [(qT)^2 - (qT)^{-2}] - \qdiff (qT - q^{-1}T^{-1}) \]

Hence $g(T) = (T^2 - T^{-2}) - \qdiff (T - T^{-1})$. Summing up, as in
the proof of Proposition \ref{polynom1}, we obtain that $\alpha_{r,m}$ equals

$$\left[ \frac{r^2(q^{2-2m}-1)q^2}{q^{-2}-1}
	- \frac{r^{-2}(q^{2m-2}-1)q^{-2}}{q^2-1} \right] $$
$$\hspace{8ex}- \qdiff \left[	\frac{r(q^{1-m}-1)q}{q^{-1}-1}-
\frac{r^{-1}(q^{m-1}-1)q^{-1}}{q-1} \right]$$

\noindent We take the ``best possible" common
factor. Then we get that
this equals
$$\frac{(q^{m-1}-1)}{q^2 r^2(q^2-1)} \Big[
	q^6r^4q^{2-2m}(q^{m-1}+1) - (q^{m-1}+1) $$
$$\hspace{8ex} - \qdiff (q^4r^3q^{1-m}(q+1) - (q+1)rq) \Big]$$

Put $q^{m-1} = T$. Then we get
\begin{align*}
& \displaystyle \frac{T-1}{q^2 r^2(q^2-1)} \big[ (T+1)((q^3 r^2
	T^{-1})^2 - 1) 
  - \qdiff qr(q+1)(q^3 r^2 T^{-1}-1) \big] \\
= & \displaystyle \frac{(T-1)}{q^2 r^2(q^2-1)} \frac{q^3 r^2 - T}{T^2}
	\left[ (T+1)(q^3 r^2 + T) - \qdiff qr(q+1) T \right] \\
= & \beta \left[ (T+1)(q^3 r^2 + T) - \qdiff qr(q+1) T \right],
	\mbox{ say.}
\end{align*}

We now show that condition (1) (of Theorem \ref{weyl}) is satisfied. If
we fix $n \in \zz$ and $r = \eps q^n$, then we want to show that there is
at most one root $m$ of the equation $\alpha_{r,m} = 0$, for this is
equivalent to condition (1), by Theorem \ref{s5}.

We know that $m \geq 2$, and $q$ is not a root of unity, hence most of
the terms in $\beta$ above are nonzero. The only term we need to consider
is $q^3 r^2 - T$. However, if $r = \eps q^n$, then this equals $q^{3+2n}
- q^{m-1}$, and for this to vanish, we need $m = 2n+4$. Clearly, this is
impossible, since we desire $\alpha_{\eps q^n,m}$ to vanish for some $2
\leq m \leq n+1$. Thus $\beta \neq 0$, so we can cancel it.

We thus need to show that if we fix $r = \eps q^n$, then there is at most
one solution of the form $T = q^{m-1}$, to the equation
\begin{equation}\label{quad}
(T+1)(q^3 r^2 + T) - \qdiff qr(q+1) T = 0
\end{equation}
where $2 \leq m \leq n+1$.

Clearly, there are no solutions when $n=0$, since $n+1 <2$. The next case
is $n=1$. The equation then becomes
\[ (T+1)(q^5 + T) = \eps T (q^4+1)(q+1). \]
We need $2 \leq m \leq 2$ to be a solution, i.e. $T = q$.

Taking $\eps = 1$,  we 
get $T^2 + q^5 = T(q + q^4)$, which holds for $T
= q,q^4$. Hence there is a unique root $T=q$, as desired.

On the other hand, if $\eps = -1$, then evaluating at $T=q$, and
cancelling $(q+1)(q^5+q)$ from both sides (since $q$ is not a root of
unity), we get $1 = -1$, a contradiction since char $\k \neq
2$. Hence, there is no root in this case.

Finally, take $n>1$. We claim, in fact, that there is no root of the
equation \eqref{quad}, of the form $T= q^{m-1}$. Simply plug in $T =
q^{m-1}$ and $r = \eps q^n$ above, and multiply both sides by $q$; we get
\[ q(q^{m-1}+1)(q^{2n+3} + q^{m-1}) = \eps (q^4+1)(q+1)q^{n + m-1}. \]
By the assumption that $q$ is transcendental, we must have the highest
degree terms on both sides to be the same. On the right hand
side, the highest degree is $4 + 1 + (n+m-1) = n+m+4$. On the left hand
side, we have $m-1 \leq n < 2n+3$, so the highest degree is $1 + (m-1) +
(2n+3) = 2n + m + 3$. These are equal only when $n=1$, so there is no
root for $n>1$.

We conclude that $V(q)$ is the unique finite dimensional simple
$A$-module (because char $\k \neq 2$). Since it has no self-extensions
(by Proposition \ref{extvan} above), every finite dimensional module is a
direct sum of copies of $V(q)$, and hence, completely reducible.
\hfill$\Diamond$

\medskip

Finally, we mention that we have similar results and (counter)examples in
the case of $H_f$ (cf. \cite{Kh}). Complete reducibility holds if and
only if every block $\calo(r)$ contains at most one finite dimensional
simple module, if and only if every $T(r)$ intersects $\Z$ in at most
four elements.

%

\section{Center}\label{scenter}

We will show in this section that the center of $A$ is trivial if $C_0
\neq 0$. Consequently, we cannot use the same approach as in \cite{BGG1}
to decompose \calo into blocks. This is why we had to follow a different
approach in Section \ref{sfinl}.

\begin{theorem} \label{thmcenter}
The center of $A$ is the set of scalars $\k$ when $C_0 \neq 0$.
\end{theorem}

\noindent The proof is in two parts. The first part is the following
proposition.

\begin{proposition}
If $a \in \mf{Z}(A)$, then $\xi(a) \in \k$, i.e. the purely CSA part of
$a$ is a scalar.
\end{proposition}

For the sake of simplicity, we denote $\mcA = \k[K,K^{-1}]$. 
Following \cite[\S 1.6]{Ja}, for $j \in \Z$, 
we define the operator $\gamma_j : \mcA \to
\mcA$ by: $\gamma_j(\varphi(K)) = \varphi(q^jK)$. Now define $\eta_j : A
\to \mcA$ by $\eta_j(a) = \gamma_j(\xi(a))$. For example, $\eta_0(C_0) =
\gamma_0(\xi(C_0)) = \xi(C_0)$.

Set $a_0 = q \cz + \cm$. We claim that $a_0 \neq 0$ if $C_0 \neq 0$. The
Casimir is $C = \xi(C) + FE$, and $C_0 = p(C) \in p(\xi(C)) + F \cdot A
\cdot E$. Hence $\xi(C_0) = p(\xi(C)) \neq 0$.

Thus, if $\cz = \xi(C_0) = \alpha_n K^n + l.o.t.$, then
$$a_0 = q \cz + \cm = \alpha_n(q + q^{-n}) K^n + l.o.t.$$

\noindent Clearly, $\alpha_n \neq 0$, so if $a_0 = 0$ we must have
$q^{-n} = -q$, or $q^{n+1} = -1$, whence $q^{2n+2} = 1$. Since $q$ is not
a root of unity, this means that $n=-1$, hence
$$a_0=2\alpha_{-1}qK^{-1}+l.o.t.\not =0,$$

\noindent which is a contradiction, since char $\k \neq 2$.

Before proving the proposition above, we need the following lemma.

\begin{lemma}
We have the following commutation relations in $A$:
\begin{enumerate}
\item If $U \in A_{q^j}$ (i.e. $KUK^{-1} = q^jU$) and $\varphi(K) \in \mcA$,
then\\ $\varphi(K) U = U \gamma_j(\varphi(K)) = U \eta_j(\varphi(K))$.

\noindent Further, when written in the PBW basis,
\item the component in $Y \cdot \mcA$ of $[X,Y^2]$ is $-Y a_0$.
\item the component in $\mcA$ of $[E,Y^2]$ is $-\cz$.
\item the component in $\mcA \cdot X$ of $[X^2,Y]$ is $-a_0 X$.
\item the component in $\mcA$ of $[F,X^2]$ is $-\cz K^{-1}$.
\end{enumerate}
\end{lemma}

\begin{proof}[Proof of the lemma]
\hfill
\begin{enumerate}
\item This is obvious.

\item We compute: $[X,Y^2] = XY^2 - Y^2X$, so the component in $Y \cdot
\mcA$ is obtained only from $XY^2$. We have $XY^2 = (qYX-C_0)Y = qY(XY) -
C_0Y = qY(qYX-C_0) - C_0Y = q^2Y^2X - qY C_0 - C_0Y$ from the defining
relations.

We need to rewrite $C_0Y$ in the PBW basis and find the component in $Y
\cdot \mcA$. Clearly, $(C_0 - \xi(C_0))Y \in A \cdot EY = A \cdot (X +
q^{-1} YE)$, and hence this contributes nothing. So the only contribution
is from $\xi(C_0)Y$, which from above equals $Y\xi(C_0)(q^{-1}K) = Y
\cm$.

In conclusion, we obtain that the desired component is $-qY \xi(C_0) -
Y \cm = -Ya_0$, as claimed.

\item Once again, we compute: $Y^2E$ can give no such component, so the
only component from $[E,Y^2]$ comes from $EY^2 = (X + q^{-1}YE)Y = XY +
q^{-1}Y(X + q^{-1} YE)$. Once again, the only contribution comes from $XY
= qYX - C_0$, and hence the component of $[E,Y^2]$ in $\mcA$ is
$-\xi(C_0) = -\cz$.

\item This is similar to above: $[X^2,Y]$ and $X^2Y$ have the same
component, which comes from $X^2Y = X(qYX-C_0) = q(XY)X - XC_0 = q^2
YX^2 - qC_0X - XC_0$. The contribution of $XC_0$ comes from
$X \xi(C_0) = \cm X$, and the contribution from $C_0X$ is $\xi(C_0)X =
\cz X$. Hence the total contribution is $-q\cz X - \cm X = - a_0 X$.

\item Finally, the component in $\mcA$ comes from $-X^2F = -X(FX -
YK^{-1}) = -XFX + (XY)K^{-1} = -XFX + qYXK^{-1} - C_0 K^{-1}$. Clearly,
only the last term has a non-zero component in $\mcA$, which is
$-\xi(C_0) K^{-1}$, as claimed.
\end{enumerate}
\end{proof}

\begin{proof}[Proof of the proposition]
Given $a \in \mf{Z}(A)$, we write $a$ as a linear combination of PBW
basis elements. Note that $KaK^{-1} = a$, whence the only basis elements
that can contribute to $a$ are of the form $F^a Y^b K^c X^d E^e$ where
$2a+b = d + 2e$.

We can write $a$ in the form
$$a = \xi(a) + Y b_1 X + Y^2 b_2 X^2 + F b_3 X^2 + F b_4 E + Y^2 b_5 E +
h.o.t.$$

\noindent Here, $h.o.t.$ denotes higher order terms in $E,X$ (i.e.
$h.o.t.$ is in the left ideal generated by $E^2, EX, X^3$) and the
$b_i$'s are Laurent polynomials in $K$.

\noindent {\bf Step 1: Obtain equations relating the coefficients $b_i$.}

We now use the fact that $a$ commutes with $X,Y,E,F$ to equate various
coefficients to zero. We have to consider six different cases.

\noindent {\bf Case 1}: The component in $Y \cdot \mcA \cdot E$, of
$[X,a]$, is zero.

Clearly, if $b \in A$, then $[X,bX] = XbX - bX^2 \in A \cdot X$, by the
PBW theorem. Similarly, $[X,bE^2] \in A E^2$, and $[X,bEX] \in A \cdot
EX$. Hence $[X,h.o.t.]$ still gives us only higher order terms. In fact,
from this analysis, we see that we only need to consider
$[X,Fb_4 E + Y^2 b_5 E]$
for the above coefficient. We have
\[ [X,Fb_4E] = [X,F] b_4 E + F[X,b_4E] = [X,F] b_4 E + F [X,b_4]
E+Fb_4[X,E], \]

\noindent and the second and third terms are clearly in $A \cdot EX$.
Hence we only need to consider the first term. The same is true for
$[X,Y^2 b_5E]$.

Hence we conclude that, to compute the above coefficient, we only need to
look at
$$[X,F]b_4 E + [X,Y^2] b_5E$$

\noindent From the lemma, the contribution is $-YK^{-1} b_4 E -Y a_0 b_5
E$. If this is to be zero, then we obtain
\begin{equation}\label{cc1}
b_4 = -K a_0 b_5
\end{equation}

\noindent {\bf Case 2}: The component in $Y \cdot \mcA \cdot X^2$ of
$[X,a]$ is zero.

Once again, by a similar analysis, we see that we only need to look at
$$Y[X,b_1]X + [X,F]b_3 X^2 + [X,Y^2]b_2 X^2$$

\noindent and the contribution is $Y[\bm - \bz]X^2 -YK^{-1} b_3 X^2-Y a_0
b_2 X^2$ from the lemma. If this is to be zero, then we obtain
\begin{equation}\label{cc2}
b_3 = K(\bm - \bz) -K a_0 b_2
\end{equation}

\noindent {\bf Case 3}: The component in $\mcA \cdot X$ of $[X,a]$ is
zero.

In this case the contribution comes from $[X,\xi(a)] + [X,Y]b_1X$. Using
the lemma, we simplify this to $(\am - \az) X -\cz b_1 X = 0$. Hence
\begin{equation}\label{cc3}
\am - \az = \cz b_1
\end{equation}

\noindent {\bf Case 4}: The component in $\mcA \cdot E$, of $[E,a]$, is
zero.

In this case we look at $([E,\xi(a)] + [E,F]b_4 + [E,Y^2]b_5)E$, which,
from the lemma above, contributes $(\amm - \az + \brac{K} b_4 - \cz
b_5)E$. If this is zero, then we get
\begin{equation}\label{cc4}
\amm - \az = -\brac{K}b_4 + \cz b_5 = (\cz + K \brac{K} a_0) b_5
\end{equation}
where the last equality follows from equation \eqref{cc1} above.

\noindent {\bf Case 5}: The component in $\mcA \cdot X^2$ of $[E,a]$ is
zero.

In this case we look at $[E,Y] b_1X + [E,Y^2]b_2 X^2 + [E,F]b_3 X^2$,
which, from the lemma above, contributes $X b_1 X - \cz b_2 X^2 +
\brac{K} b_3$. If the contribution from this is zero, then we get
\begin{equation}\label{cc5}
\cz b_2 = \bm + \brac{K} b_3
\end{equation}

\noindent {\bf Case 6}: The component in $F \cdot \mcA \cdot X$ of
$[Y,a]$ is zero.

In this case the contribution comes from $-Fb_3[X^2,Y] -F b_4[E,Y]$.
Using the lemma, we simplify this to $Fa_0 b_3 X - Fb_4 X = 0$. Hence
\begin{equation}\label{cc6}
b_4 = a_0 b_3
\end{equation}\hfill

\noindent {\bf Step 2: Solve the above system for the $b_i$'s.}

We now use these equations. From equations \eqref{cc1} and \eqref{cc6},
we get that $a_0(b_3 + Kb_5) = 0$. We proved at the beginning of this
section that $a_0 \neq 0$. Hence $b_3 = -Kb_5$.

Multiplying equation \eqref{cc2} by $\cz$, and using equation
\eqref{cc5}, we get
$$\cz b_3 = \cz K(\bm - \bz) - Ka_0(\bm+\brac{K} b_3)$$
so that
$$(\cz + K \brac{K} a_0) b_3 = -K[(a_0-\cz) \bm + \cz \bz]$$
and this equals $-K (\cz + K\brac{K} a_0) b_5$ because $b_3 = -Kb_5$.
Using equation \eqref{cc4}, and dividing by $-K$, we finally get
\begin{eqnarray*}
\amm - \az & = & (a_0 - \cz) \bm + \cz \bz\\
           & = & [(q-1) \cz + \cm] \bm + \cz \bz\\
           & = & (q-1) \cz \bm + (\cm \bm + \cz \bz)
\end{eqnarray*}

\noindent Thus we finally get, using \eqref{cc3},
$$\amm - \az = (q-1) \cz \bm + (\amm - \am) + (\am - \az)$$
so that
$$(q-1) \cz \bm = 0$$
The above holds in $\mcA$. Since $(q-1) \cz = (q-1) \xi(C_0) \neq 0$ by
assumption, $\bm = 0$. Finally, applying $\eta_{-1}$ to equation
\eqref{cc3}, we get that $\amm = \am$. But if $\xi(a) = \sum_i \alpha_i
K^i$, then this gives $\alpha_i q^{-i} = \alpha_i q^{-2i}$ for all $i$.
Since $q$ is not a root of unity, the only nonzero coefficient is
$\alpha_0$ and $\xi(a) = \alpha_0$ is indeed a scalar, as claimed.
\end{proof}

To complete the proof that the center is trivial, we use the PBW form of
the basis. The lemma below says that for any ``purely non-CSA" element
$\beta \neq 0$, we can find $w_r \in Z(r)$ (for ``most" $r \neq \pm q^n$)
so that $\beta w_r \neq 0$ in $Z(r)$. In fact, we explicitly produce such
a $w_r$.

Suppose we are given $\beta \in A$ so that $\xi(\beta) = 0$, and $\beta
\neq 0$. We can write $\beta$ in the PBW form $\beta = \sum_i \beta_i
p_i(K) X^{d_i} E^{e_i}$. Here, $\beta_i \in \k[Y,F]$ and $p_i$'s are
Laurent polynomials in one variable. Choose $i$ so that $e = e_i$ is the
least among all $e$'s, and among all $j$'s with $e_j = e$, the least
value of $d_j$ is $d = d_i$. Without loss of generality, 
we may assume $i=0$.

\begin{lemma}
There exists a finite set $T \subset \k$ with $0\in T$ such that if $r
\neq \pm q^n$, $r \notin T$ and if $w_r = F^e v_{q^{-d}r}$, then $\beta
w_r \in \k^\times \beta_0 v_r$.
\end{lemma}

\begin{proof}
We work in the Verma module $Z(r)$, where $r \neq \pm q^n$ for any $n
\geq 0$. We define $w_r = F^e v_{q^{-d}r}$ and compute $X^{d_i} E^{e_i}
w_r$.

Since $v_{q^{-d}r}$ is annihilated by $E$, it generates a
\uqsl-Verma module $Z_C(q^{-d}r)$, and by \uqsl-theory we
observe that $E^e F^e v_{q^{-d}r}$ is a nonzero scalar multiple of
$v_{q^{-d}r}$ (by \cite[Proposition 2.5]{Ja} the Verma module
is simple, so the only vector killed by $E$ is $v_{q^{-d}r}$).

Next, using equation \eqref{Rn}, an easy induction argument shows that
\begin{equation}\label{centerlemma}
X^d v_{q^{-d}r} = (-1)^d \prod_{i=1}^d
\frac{\alpha_{r,d+2-i}}{\brac{q^{1+i-d}r}}\quad v_r
\end{equation}

For each fixed $i$, the expression $\alpha_{r,d+2-i}$ is a nonzero
(Laurent) polynomial in $r$, hence it has a finite set of roots. We now
define the finite set $T$ of ``bad points". Recall that we wrote $\beta =
\sum_i \beta_i p_i(K) X^{d_i} E^{e_i}$. Define $T$ to be the union of the
(finite) set of roots of $p_0$, 0, and the (finite) set of roots $r$ of
all the $\alpha_{r,d+2-i}$ for $1\le i\le d$.

Finally, we compute $X^{d_i} E^{e_i} w_r$. There are two cases:

(a) $e_i > e$, in which case  $E^{e_i} w_r = E^{e_i - e - 1} (E^{e+1} F^e
v_{q^{-d}r}) = 0$ by \uqsl-theory; or

(b) $e_i = e$ ($i=0$), in which case $X^{d_i} E^e w_r = X^{d_i - d} (X^d
E^e F^e v_{q^{-d}r})$. From above, if $r \notin T$, then this is
$X^{d_i-d} c v_r$ for some nonzero scalar $c$. Thus, we get a nonzero
vector if and only if $d_i = d$ since $v_r$ is maximal.

Thus, $\beta w_r = c \beta_0 p_0(K) v_r = c \beta_0 p_0(r) v_r$. Hence
$\beta w_r = (c p_0(r)) (\beta_0 v_r)$ and $cp_0(r) \neq 0$ for all $r
\notin T,\ r \neq \pm q^n$.
\end{proof}

\begin{proof}[Proof of Theorem \ref{thmcenter}]\hfill

Suppose $a = \xi(a) + \beta \in \mf{Z}(A),\beta\not\in \mcA$ and $\beta
\neq 0$. Let us look at how $a$ acts on $w_r = F^e p_{d,r}(Y,F) v_r$ (as
above), with $r \notin T$ and $r\not=\pm q^n$. We know $\beta w_r = f(r)
\beta_0 v_r$, $f(r)\in\k^{\times}$. Now, $a(F^e p_{d,r}) = (F^e
p_{d,r})a$, since $a$ is central. Thus, $a w_r = F^e p_{d,r}(Y,F) av_r$,
i.e. $\xi(a)w_r + \beta w_r = F^e p_{d,r}(Y,F) \xi(a) v_r + F^e
p_{d,r}(Y,F) \beta v_r = \xi(a)(r) w_r + 0 = \xi(a)(r)w_r$.

\noindent Thus, $f(r) \beta_0v_r = (\xi(a)(r) - \xi(a)(q^{-n}r))w_r$ for
some $n$, i.e.
$$(\xi(a)(r) - \xi(a)(q^{-n}r)) F^e p_{d,r}(Y,F) = f(r)\beta_0, \mbox{
for all } r\not\in T,r\not=\pm q^{n}.$$
But from the above proposition, $\xi(a)$ is a constant, so $\beta_0 = 0$
because $f(r)\not=0 $. This contradicts our assumption that $\beta\not=
0$. Therefore, $\beta=0$ and we conclude that $a=\xi(a)\in\k^{\times}$ so
that the center is trivial.
\end{proof}

\section{Counterexamples}\label{scex}

We provide counterexamples for two questions:

\begin{enumerate}
\item Is every Verma module $Z(r)$ a direct sum of \uqsl-Verma modules
$$Z_C(r) \oplus Z_C(q^{-1}r) \oplus \dots ?$$
\item If $\alpha_{r,n+1} = 0$, is it true that $Z(q^{-n}r)
\hookrightarrow Z(r)$ ?
\end{enumerate}

The answers to both questions are: no.
\begin{enumerate}
\item The structure equations guarantee, for $r = \eps q^n$, that
$v_{\eps q^{-1}}$ can be defined, and is \uqsl-maximal. However, if
$Z(r)$ is to decompose into a direct sum of $Z_C(r')$'s (as above), then
we need a monic polynomial $h(Y,F) = Y^{n+2} + l.o.t.$, so that $v_{\eps
q^{-2}} = h(Y,F)v_r$ is \uqsl-maximal.

Now, $EYv_{\eps q^{-1}} = Xv_{\eps q^{-1}} = -\alpha_{r,n+2} v_\eps$, by
equation \eqref{Rn}. By \uqsl-theory, $EF^lv_{\eps q^{2l}} \in \k
F^{l-1} v_{\eps q^{2l}}$ for each $l>0$. Thus, if there exists a
\uqsl-maximal vector, it has to be a linear combination of $Yv_{\eps
q^{-1}}$ and $Fv_\eps$. However, $EFv_\eps = 0$, so the only way $Y^{n+2}
+ FY^n + l.o.t.$ is \uqsl-maximal, is if $\alpha_{r,n+2} = 0$. By
definition of $\alpha$, this holds if and only if $\alpha_{r,n+3} = 0$.

We conclude that $Z(\eps q^n)$ has a \uqsl-Verma component
$Z_C(\eps q^{-2})$ only if $\alpha_{r,n+3} = 0$. Hence (1) fails in
general.

\item This requires some calculations. By definition, we see that
$\alpha_{\eps, 4} = 0$. We now show that $Z(\eps)$ does not always have a
Verma submodule $Z(\eps q^{-3})$.

By Proposition \ref{s8}, if there exists a maximal vector of weight
$\eps q^{-3}$, then (up to scalars) it {\it must} be $v' = v_{\eps q^{-3}}
= (Y^3 - b FY) v_\eps$, where
$$b = \eps((q+q^{-1})c_{0,\eps} + c_{0, \eps q^{-1}})$$

We now calculate what happens when this vector is also killed by $X$.
From the proof of Proposition \ref{s8}, we know that $Xv' = b'F v_\eps$,
because the coefficient of $Y^2v_\eps$ was made to equal zero. We now
show that $b'$ is not always zero.

Clearly,
$$XFYv_\eps = (FX-YK^{-1})Yv_\eps = F(XYv_\eps) - \eps q^{-1} Y^2
v_\eps$$

\noindent and
$$F(XYv_\eps) - \eps q^{-1}Y^2 v_\eps = -(Fc_{0,\eps} + \eps q^{-1} Y^2)
v_\eps.$$

But \begin{align*}
XY^3 v_\eps =& (qYXY^2 - C_0 Y^2)v_\eps \\
=& (q^2Y^2 XY - q Y C_0 Y -
C_0 Y^2) v_\eps \\
=& - q^2 c_{0,\eps} Y^2 v_\eps - q Y c_{0,\eps q^{-1}} Y
v_\eps - C_0 Y^2 v_\eps.
\end{align*}
Hence we only need to look at $C_0 Y^2 v_\eps$, to find the coefficient
of $Fv_\eps$.

The basic calculation is this: $E Y^2 v_\eps = XYv_\eps = -c_{0,\eps}
v_\eps$. Hence,
$$C Y^2 v_\eps = -c_{0,\eps} Fv_\eps + c_{\eps q^{-2}} Y^2
v_\eps = -c_{0,\eps} Fv_\eps + c_{\eps} Y^2 v_\eps,$$ 

\noindent by definition of $c_r$. An easy induction argument now shows
that 
\begin{align*}
C_0 Y^2 v_\eps =& p(C) Y^2 v_\eps \\
=& \frac{-c_{0,\eps}}{c_\eps -
c_{0,\eps}} (p(c_\eps) - p(c_{0,\eps}))F v_\eps + p(c_\eps) Y^2 
v_\eps\\
=& \frac{-c_{0,\eps}}{c_\eps - c_{0,\eps}}(c_{0,\eps} - p(c_{0,\eps}))
Fv_\eps + c_{0,\eps} Y^2 v_\eps \\
=& \frac{-c_{0,\eps}}{c_\eps - c_{0,\eps}}(c_{0,\eps} - p(c_{0,\eps}))
Fv_\eps + c_{0,\eps} Y^2 v_\eps \\
=& -a Fv_\eps + c_{0,\eps} Y^2 v_\eps, \mbox{ say.}
\end{align*}

Hence, we conclude that the coefficient of $Fv_\eps$ in $Xv' = X(Y^3 -
bFY) v_\eps$ is $b c_{0,\eps} + a$, and this should be zero if $v'$ is
maximal. Simplifying, we get
\begin{equation}\label{suff}
c_{0,\eps}(c_\eps - c_{0,\eps})\eps((q+q^{-1})c_{0,\eps} + c_{0,\eps
q^{-1}}) + c_{0,\eps}(c_{0,\eps} - p(c_{0,\eps})) = 0
\end{equation}

But this is not always satisfied: take $p(T) = \beta T$ for some $\beta
\in \k,\ \beta \neq 0,1$. Then the above condition reduces to
$$\frac{(q+q^{-1})^2 + 2}{(q-q^{-1})^2} + 1 = 0$$

\noindent which simplifies to $2q^6 = 2$. However, since char $\k \neq 2$
and $q$ is not a root of unity, this is not true. So at least for some
$p(C)$, this condition is false.
\end{enumerate}

\section{Classical limit}\label{sclaslimit}


The algebra $A$ specializes to the symplectic oscillator algebra $H_f$ of
\cite{Kh} as $q\to 1$; this is what we formalize in this section. 
Let $k$ be a field of characteristic 0.
Let $\k = \kk(q)$ be the field of rational functions on $\kk$
and let $R \subset \k$ be the $\kk$-subalgebra of rational functions
regular at the point $q=1$.
Recall from \cite{Kh} that
$$\Delta_0 := 1+f(\Delta), \quad \Delta := (FE+H/2+H^2/4)/2,$$

\noindent where $f\in k[t]$. $H_f$ is the $\kk$-algebra with generators
$X,Y,E,F,H$ with relations: $\langle E,F,H\rangle$ generate \usl,
$[E,X]=[F,Y]=0,[E,Y]=X,[F,X]=Y,[H,X]=X,[H,Y]=-Y$ and $[Y,X]=1+f(\Delta)$.

We write $\Delta_0$ as $$\Delta_0 = f_0(FE+(H+1)^2/4)$$ for some $f_0$, a
polynomial in one variable with coefficients in $\kk$. We will explain
how $H_f$ is the limit of $A$ as $q \to 1$.

Our algebra $A$ is fixed, and in particular, so is the polynomial $p$.
Write $C_0$ as
$$ C_0 = f_0\Big(FE+\frac{Kq+K^{-1}q^{-1}-2}{(q-q^{-1})^2}\Big)$$

\noindent for some polynomial $f_0$. The coefficients of $f_0$ are in
$\kk$, but the limiting process works so long as they are in $R$. We will
follow the approach in \cite{HK} and use the notation on pp. 48; in
particular,
$$ (K^m;n)_q := \frac{K^mq^n-1}{q-1} ,\quad m,n\in \Z.$$

We define $A_R$ to be the $R$-subalgebra of $A$ generated by the elements
$X, Y, E, F, K, K^{-1}$, $(K;0)_q$, and set
\begin{equation} \label{AA1}
A_1 := (R/(q-1)R)\otimes_{R}A_R=A_R/(q-1)A_R.
\end{equation}
The elements $(K^m;n)_q$ are
all in $A_R$. This happens in the case $n=0$, because
\[ K^m(K;0)_q = (K^{m+1};0)_q - (K^m;0)_q, \]
so by induction it follows that $(K^m;0)_q\in A_R$. For general $n$, we
now conclude that
\[ (K^m;n)_q = K^m(1;n)_q + (K^m;0)_q \in A_R. \]\hfill

\begin{proposition}
The algebra $A_1$ (defined in (\ref{AA1}))
is isomorphic to $H_f$.
\end{proposition}

\begin{proof}
Denote by $\oX,\oY,\oE,\oF,\oK^m$ ($m\in\Z$) the images of $X,Y,E,F,K^m$
under $A_R \to A_1$. Then the image of $K^m-1$ equals the image of $(q-1)
(K^m;0)_q \in A_R$. Thus $K^m-1 \mapsto 0$ in $A_1$, so $K^m \mapsto 1$
under $A_R\to A_1$, for all $m \in \Z$.

Define the element $\oH$ in $A_1$ to be the image of $(K;0)_q$ under the
projection $A_R \to A_1$. The element $C_0$ is in $A_R$: this follows
from the observation that
\begin{equation} \label{eqkq} 
\frac{Kq+K^{-1}q^{-1}-2}{(q-q^{-1})^2} = \frac{K^{-1}q
(K;1)_q^2}{(q+1)^2}
\end{equation}

\noindent is in $A_R$, which in turn is a consequence of \cite[Lemma
3.3.2]{HK}. The image of $\frac{K^{-1}q (K;1)_q^2}{(q+1)^2}$ under the
projection $A_R\to A_1$ is $(\oH+1)^2/4$:
$$ \frac{K^{-1}q (K;1)_q^2}{(q+1)^2} = \frac{K^{-1}q}{(q+1)^2}(q(K;0)_q 
+1)^2 $$

\noindent and we know that $K\to 1$ and $(K;0)_q\to\oH$.

Therefore, because of our choice of $f$, $\oX$ and $\oY$ satisfy the
relation
\[ \oY\oX-\oX\oY=f_0(\oF\oE+(\oH+1)^2/4).\]

It is clear that, in $A_1$, we have the relations
$\oE\oX=\oX\oE, \oE\oY-\oY\oE=\oX,\oF\oX-\oX\oF=\oY$ and $\oF\oY=\oY\oF$.
Therefore, we have an epimorphism $H_f\to A_1$. These two rings have a
filtration where $\ddeg(X)=\ddeg(Y)=1,\ \ddeg(E)=\ddeg(F)=\ddeg(H)=0$ and
similarly with $\oX,\oY,\oE,\oF,\oH$.

We can identify $\oE$ with $E$, $\oF$ with $F$ and $\oH$ with $H$,
because we know from \cite{HK} that $\oE,\oF,\oH$ generate a subalgebra
isomorphic to \usl. 

We can view the map $H_f\to A_1$ as a map of
\usl-modules. Now, $gr(H_f)= k[X,Y]\rtimes \usl$.
Also, $gr(A_1)= k[X,Y]\rtimes \usl$ since
$gr(A_R)=R[X,Y]\rtimes \uqsl$.
The associated graded map 
$gr(H_f)\to gr(A_1)$
is the identity map from
$k[X,Y]\rtimes \usl$ to
$k[X,Y]\rtimes \usl$.
Hence, $H_f$ is isomorphic to $A_1$.
\end{proof}
 
Let $r=\pm q^n$ where $n\in\Z$ and let $V$ be a standard cyclic
$A$-module with highest weight $r$ and highest weight vector $v_r$. We
define the $R$-form of $V$ to be the $A_R$-module $V_R := A_R \cdot v_r$.
Set $V^1 := R/(q-1) \otimes_R V_R,$ so $V^1$ is an $A_1$-module.

\begin{proposition}
$V_1$ is an $H_f$ standard cyclic module with highest weight $n$ and
highest weight vector $v_r$. Furthermore, if $V$ is a Verma module, then
so is $V_1$.
\end{proposition}                 

\begin{proof}
This is clear from the previous proposition.
\end{proof}

\begin{acknowledgement}
We thank Pavel Etingof, Victor Ginzburg, Anna Lachowska, and Alexei
Oblomkov for stimulating discussions. We also thank the referee
for his suggestions and comments.
The work of W.L. Gan was partially supported by
the NSF grant DMS-0401509.
\end{acknowledgement}

\end{document}